\documentclass[11pt]{article}
\usepackage{graphicx}
\usepackage{amssymb}
\usepackage{epstopdf}
\newtheorem{theorem}{Theorem}[section]
\newtheorem{corollary}[theorem]{Corollary}
\newtheorem{proposition}{Proposition}[section]
\newtheorem{definition}{Definition}[section]
\newtheorem{lemma}{Lemma}[section]

\title{ A Probabilistic Representation of Solutions of the 
Incompressible Navier-Stokes Equations in
\({\bf R}^3\)\footnote{
words and 
phrases: tree-indexed Markov process, 
branching random walk,
incompressible
Navier-Stokes equations.} }
\author{Mina Ossiander\footnote{This research was partially supported by NSF FRG  collaborative award DMS-0073958.} \\ Department of Mathematics \\ Oregon State University }
\begin{document}
\maketitle{}

\maketitle{}

\begin{abstract}  A new probabilistic representation is presented for solutions of the incompressible Navier-Stokes equations in \({\bf R}^3\) with given forcing and initial velocity.  This representation expresses solutions as scaled conditional expectations of functionals of a Markov process indexed by the nodes of a binary tree.  It gives existence and uniqueness of weak solutions for all time under relatively simple conditions on the forcing and initial data.  These conditions involve comparison of the forcing and initial data with majorizing kernels.
\end{abstract}

\section{Introduction.}

The 3-dimensional Navier-Stokes equations governing the velocity of incompressible fluids with viscosity \(\nu\) and initial velocity \(u_0\) are given by
\begin{eqnarray}\label{ns}
{\partial u \over \partial t} + u \cdot \nabla u & = & \nu  \Delta u - \nabla p + g \nonumber \\
\nabla \cdot u & = & 0 \\
u_0(x) & = & u(x,0) \nonumber
\end{eqnarray}
Here \(u:  {\bf{R}}^3 \times [0, \infty )  \to {\bf{R}}^3 \) is the velocity vector, \(p:  {\bf{R}}^3 \times [0, \infty )  \to {\bf{R}} \) is the pressure and \(g:  {\bf{R}}^3 \times [0, \infty )  \to {\bf{R}}^3 \) is the external forcing.   Although these equations have been studied extensively over the last century important open questions remain concerning existence and uniqueness of smooth solutions.
Extensive background and development of results can be found in Temam \cite{T79}, Doering and Gibbon \cite{DG95}, Foias, Manley, Rosa and Temam \cite{FMRT01}, and Lemari\'{e}-Rieusset \cite{L-R02}.  Ladyzhenskaya \cite{L03} provides a recent interesting and accessible review of results.

There are two classes of existence and uniqueness results for solutions to the Navier-Stokes equations; existence and uniqueness for all time if the forcing and initial data are small enough in some specified function space, or existence and uniqueness for a short time for larger initial data.  This paper falls into the former camp in giving a probabilistic representation in physical space of the solutions of (\ref{ns}) when a smallness condition is satisfied by the given forcing and initial velocity.  This results in a method of demonstrating existence and uniqueness of weak solutions (in the sense of Lemari\'{e}-Rieusset \cite{L-R02}) for all time under straightforward bounds on forcing and initial data.  The focus is strictly on the Navier-Stokes equations in \({\bf R}^3\).  Indeed the method developed relies heavily on an integral formulation of the Navier-Stokes equations in \({\bf R}^3\) incorporating Leray's projection onto the subspace of divergence free vector fields.  

The stochastic representation given here has parallels with the representations given by McKean \cite{McK75} for solutions of the Kolmogorov-Petrovskii-Piskunov equation and Le Jan and Sznitman \cite{LJS97} for solutions of the Fourier transformed Navier-Stokes equations.  In the seminal paper \cite{McK75} McKean gives the solution of the KPP equation
\[
\frac{\partial u}{\partial t} = \frac{1}{2} \, \frac{\partial^2 u}{\partial x^2} + u^2 -u 
\]
with initial data \(u(x,0)=f(x) \), \( 0 \leq f(x) \leq 1 \) for all \(x \in {\bf R}\), as
\[
u(x,t) = E_x \prod_{i=1}^{N(t)} f(W^{(i)}(t))
\]
where \({\bf W} = \{W^{(i)} : i \geq 1 \}\) is a branching Brownian motion in \(\bf R\) with initial value \(x\).  Roughly here one can think of the Laplacian term \(\frac{1}{2} \, \frac{\partial^2 u}{\partial x^2}\) as producing the Brownian motion, the \(u^2\) term yielding the binary branching, and the \(-u\) term resulting in the exponential waiting time between branches.  The existence of the solution is a consequence of the bound imposed on the initial data.  

LeJan and Sznitman in \cite{LJS97} give the solution of the Fourier transformed Navier-Stokes equations
in \({\bf R}^3\)
\begin{eqnarray}\label{eqFNS}
 \hat{u}( \xi ,t) & = & e^{- \nu | \xi |^2t } \hat{u}_0( \xi ) \\
 & + &  \int_{s=0}^t (2 \pi)^{-3/2} e^{- \nu | \xi |^2 s}
 \Big\{-i {\bf{P}}_\xi  \int_{{\bf{R}}^3} \hat{u}( \eta, t-s) \, \xi \cdot \hat{u} ( \xi - \eta , t-s)d \eta
 + {\bf{P}}_\xi \hat{g}( \xi ,t-s) \Big\} ds \nonumber
 \end{eqnarray}
in the form of a scaled expectation 
 \[
 \hat{u}(\xi,t) = H( \xi ) E_{\xi,t} {\bf M} 
 \]
 if both
 \begin{equation}\label{eqFscale}
 \frac{| \hat{u}_0 ( \xi)|}{H(\xi)}   \quad \hbox{and} \quad 
 \frac{|\hat{g}(\xi,t)|}{|\xi|^2 H(\xi)}
 \end{equation}
 are uniformly small enough.   Here \(\bf M\) is a multiplicative functional of a stochastic cascade rooted at \(\xi\) at time \(t\) with inputs derived from initial data and forcing scaled as in (\ref{eqFscale}) above.  The scalar function \(H:{\bf R}^3 \to (0,\infty]\) is given by either \(H(\xi)=|\xi|^{-2}\) or \(H(\xi)=\alpha |\xi|^{-1} e^{-\alpha |\xi|}\) for some \(\alpha > 0\).  The vector projection \({\bf P}_\xi\) onto the plane perpendicular to \(\xi\) has eliminated the pressure \(p\).  The term \(\hat{u}( \eta, t-s) \, \xi \cdot \hat{u} ( \xi - \eta , t-s)\) yields binary branching in space according to the normalized Markovian kernel \(\frac{H(\eta)H(\xi-\eta)}{|\xi|H(\xi)}\)  and the coefficients \(e^{- \nu | \xi |^2t } \) result in exponential waiting times between branches.  This general framework is elaborated and extended via the introduction of Fourier multiplier majorizing kernels in Bhattacharya et al \cite{BCDGOOTW03}  and \cite{BCGOOTW04}.  Again the existence of solutions is a consequence of the functional  bounds imposed on the forcing and initial data.
 
The goal of this article is to construct a binary branching process with jumps that corresponds naturally to a formulation of solutions in physical space to the Navier-Stokes equations.  It is not too surprising that this stochastic representation is 
more complicated than representations of solutions to the KPP equation or the Fourier transformed Navier-Stokes equations.  Majorizing kernel pairs \((h, \tilde{h})\) are introduced to be used as scaling multipliers.  As such they both dominate the forcing and velocity and allow the construction of transition densities.  Additionally, combined with the Laplacian term \(\nu  \Delta u\), their presence can give rise to a \(h\)-Brownian motion.  The binary nature of the branching  derives from the term \(u \cdot \nabla u\).  In this representation, the waiting time between branches is related to first passage times for Brownian motion rather than being exponential in nature.  This somewhat complicated relationship is in part due to the incompressibility constraint.  This leads to the construction of a Markov process indexed by the binary tree rather than a Markov process in time.  Analogous to the work of LeJan and Sznitman \cite{LJS97} and Bhattacharya et al \cite{BCDGOOTW03}, \cite{BCGOOTW04}, the existence of solutions is a consequence of the bounds imposed on the forcing and initial data via majorizing kernel domination.  The following theorem is representative. 
\begin{theorem}\label{thm1} Let \(h:{\bf R}^3 \to (0,\infty]\) and \(\tilde{h}:{\bf R}^3 \to [0, \infty]\) with \(h\) locally square integrable and \(h\), \(\tilde{h}\) jointly satisfying 
\begin{equation}\label{eqconth1}
\int_{{\bf R}^3} h^2(x-y) |y|^{-2}dy \leq h(x) \quad \mbox{and} \quad
\int_{{\bf R}^3} \tilde{h}(x-y) |y|^{-1}dy \leq h(x) \quad \mbox{for all} \quad x \in {\bf R}^3.
\end{equation}
If for all \(x \in {\bf R}^3\) and \(t>0\),
\[
(4 \pi \nu t)^{-3/2}| \int_{{\bf R}^3} u_0(x-y)e^{-|y|^2/4 \nu t}dy |< \pi \nu h(x)/11 \quad \mbox{and} \quad
|g(x,t)| < (\pi \nu)^2 \tilde{h}(x)/11,
\]
then there exists a collection of probability measures \(\{P_x : x \in {\bf R}^3\}\) defined on a common measurable space \((\Omega, {\cal F})\), a measurable function \(\Xi : (0,\infty) \times \Omega \to {\bf R}^3\),
 and a weak solution \(u(x,t)\) to the Navier-Stokes equations (\ref{ns}) with 
\[ P_x( \{ \omega : |\Xi(t, \omega)| < 2 \pi \nu /11 \quad \mbox{for all} \; \; t>0 \})=  1 \quad \mbox{for all} 
\quad x \in {\bf R}^3 
\]
and
\[
u(x,t)= h(x) \int_\Omega \Xi(t,\omega) dP_x(\omega) \quad \mbox{for all} \quad x \in {\bf R}^3, t>0.
\]
Furthermore, the solution \(u\) is unique in the class 
\(\{ v \in ({\cal S}'( {\bf R}^3 \times (0, \infty)))^3 : |v(x,t) | < 2 \pi \nu h(x) /11\) for all \(x \in {\bf R}^3, t>0 \}.\)
\end{theorem}
The probabilistic representation gives the scaled solution of the Navier-Stokes equations as the expectation of \(\{\Xi(t): t \geq 0 \} \), which is defined in Sections 4 and 5 as a function of a 
Markov process indexed by the nodes of a binary tree.  
The scaling, \(h\), is a non-negative majorizing kernel that controls the magnitude of the initial data \(u_0\).  The kernel \(h\) is paired with a second kernel \(\tilde{h}\) that controls the magnitude of the forcing \(g\).   The conditions stipulated by (\ref{eqconth1}) allow the definition of stochastic transition probabilities for an underlying branching Markov random walk while controlling the magnitude of the function \(\Xi\).  For example, two majorizing kernel pairs satisfying (\ref{eqconth1}) are \((h(x),\tilde{h}(x))= (\pi^{-3}|x|^{-1},(2 \pi^2|x|(1+|x|))^{-2})\) and \((h(x),\tilde{h}(x))= (\pi^{-3}(1+|x|)^{-1}, (4 \pi^4|x|(1+|x|)^3)^{-1})\).  
This representation also gives a connection between Brownian motion in \({\bf R}^3\) and solutions to  the Navier-Stokes equations.

The results of this paper are largely self-contained.  The organization is as follows.  Section 2 gives a definition of the solution space involved along with the derivation of an integral formulation of the incompressible Navier-Stokes equations incorporating Leray's projection.  The integral formulation is given explicitly in terms of the transition density of Brownian motion in \({\bf R}^3\).  Majorizing kernel pairs are defined and some important properties of classes of majorizing kernels are derived in 
Section 3.  Section 4 gives a probabilistic representation of the scaled weak solution to (\ref{ns}) in the case of excessive majorizing kernels.  This representation involves Doob's \(h\)-Brownian motion in \({\bf R}^3\).  This section also contains existence and uniqueness results in the case that the initial data \(u_0\) is dominated appropriately by an excessive kernel \(h\).  An alternate probabilistic representation for general classes of majorizing kernels is given in Section 5.  The existence and uniqueness results given in this section are consequently more general as well.  The proof of Theorem \ref{thm1} is included in Section 5.  Section 6 contains a few concluding remarks.

\section{An integral reformulation of the incompressible Navier-Stokes equations.} 

The delineation of the definition of solutions and the function spaces in which solutions exist is as follows.  
Throughout elements of \({\bf{R}}^3\) are treated as \(3 \times 1\) vectors.  The usual vector dot product and Euclidean norm \(| \cdot |\) are used.  \({\cal D}({\bf R}^3) = {\cal C}^\infty_0({\bf R}^3)\) and \({\cal D}({\bf R}^3 \times (0, \infty))
 = {\cal C}^\infty_0( {\bf R}^3 \times (0, \infty))\) where \({\cal C}^\infty_0\) is the space of \({\cal C}^\infty\) functions with compact support.  \({\cal S}'( {\bf R}^3 \times (0, \infty))\) is the space of tempered
distributions on \({\bf R}^3 \times (0, \infty)\).   
\begin{definition}  A weak solution of the Navier-Stokes equations with initial velocity \(u(x,0)=u_0(x)\) is a vector field \(u(x,t) \in
({\cal S}'( {\bf R}^3 \times (0, \infty)))^3\) satisfying the following:

\parskip = 0.05 in
\noindent
(i.)  u is locally square integrable on \({\bf R}^3 \times (0, \infty)\),

\noindent
(ii.)  \(\nabla \cdot u =0\), and

\noindent
(iii.)  there exists \(p \in {\cal S}'( {\bf R}^3 \times (0, \infty))\) with 
\({\partial u \over \partial t} + u \cdot \nabla u = \nu  \Delta u - \nabla p + g\).
\end{definition}

After incorporating incompressibility via the Leray projection \(\cal P\)
and applying Duhamel's principle, 
 the pressure term can be dropped and the Navier-Stokes equations written in integral form as 
\begin{eqnarray}\label{ins}
u & = &  e^{ \nu t \Delta } u_0 - \int_0^t e^{\nu (t-s) \Delta } \,  {\cal P} \nabla \cdot ( u \otimes u) (s) ds 
+ \int_0^t e^{\nu (t-s) \Delta} {\cal P} g(s) ds \\
\nabla \cdot u_0 & = & 0 \nonumber
\end{eqnarray}
where the Leray projection \(\cal P\) onto the space of divergence free vector fields of a vector-valued function \(v\) is defined as
\[
{\cal P} v = v - \nabla \cdot {\Delta}^{-1}( \nabla v).
\]
See Foias et al \cite{FMRT01} or Lemari\'{e}-Rieusset \cite{L-R02} for discussion.
The integral formulation above is 
expressed explicitly in Proposition \ref{propins} that follows.
\begin{definition}  A function \(h:{\bf R}^3 \to {\bf R}\) is uniformly locally square integrable if for all \(\varphi \in {\cal D}({\bf R}^3)\), 
\[
\sup_{x \in {\bf R}^3} \int_{{\bf R}^3} |\varphi(x-y) \,  h(y)|^2 dy < \infty.
\]
\end{definition}
\begin{definition}  A weak solution \(u\) of the Navier-Stokes equations on \({\bf R}^3 \times (0, \infty)\) is uniformly locally square integrable if for all \(\varphi \in {\cal D}({\bf R}^3 \times (0,\infty))\) 
\[
\sup_{x \in {\bf R}^3}  \int_{{\bf R}^3} \int_{(0,\infty)} |\varphi(x-y,t) \, u(y,t)|^2 dy dt < \infty.
\]
\end{definition}
The following equivalence theorem can be found in Lemari\'{e}-Rieusset \cite{L-R02}.
\begin{theorem}  If \(u\) is uniformly locally square integrable, then \(u\) is a weak solution of (\ref{ns}) if and only if \(u\) is a solution of (\ref{ins}).
\end{theorem}
The following definitions are useful in stating Propostion \ref{propins}.  Let
\[
K(y,t) = (2 \pi t)^{-3/2} e^{-|y|^2/2t}, \quad t>0, \quad  y \in {\bf{R}}^3
\]
denote the transition 
density of a 3-dimensional Brownian motion starting at the origin at time \(t=0\).  
For \(y \in {\bf{R}}^3\) with \(|y| >0\), let  \(e_y =y/|y|\) denote the unit vector in the direction of \(y\) and let  \({\bf{P}}_y\) denote the 3 by 3 matrix that projects vectors onto the space perpendicular to 
\(y\); so that the  entries of \(  {\bf{P}}_y \) are  
\[
({\bf{P}}_y)_{i,j} =\delta_{i,j} -(e_y)_i (e_y)_j.
\]
Define the bilinear forms \({\bf b}_1\) and \({\bf b}_2\) on \({\bf{R}}^3 \setminus \{(0,0,0)\}  \times ({\bf{R}}^3)^2\) via
\[
{\bf b}_1 (y;u,v) =(u \cdot e_y) {\bf{P}}_y v +(v \cdot e_y) {\bf{P}}_y u
\]
and
\[
{\bf b}_2 (y;u,v) = {\bf b}_1(y;u,v) + u \cdot(I-3e_y e^t_y)v \, e_y.
\]
\begin{proposition}\label{propins}
If \(u:{\bf R}^3 \times (0, \infty) \to {\bf R}^3\) is locally square integrable and satisfies
\begin{eqnarray}\label{eqexpins}
u(x,t) & = & \int_{{\bf R}^3} u_0(x-y)K(y,2 \nu t) dy \nonumber \\
         & + & \int_{s=0}^t \int_{{\bf R}^3}  \Bigg\{ \frac{|z|}{4 \nu s} K(z,2 \nu s) \;
 {\bf b}_1(z;u(x-z,t-s),u(x-z,t-s))  \\
 & + & \Big( \frac{1}{|z|} K(z,2 \nu s) - \frac{3}{4 \pi |z|^4} \int_{\{y:|y| \leq |z| \}}
K(y, 2 \nu s) dy \Big) \; {\bf b}_2(z;u(x-z,t-s),u(x-z,t-s)) \nonumber \\
 & + & \Big( K(z, 2 \nu s) {\bf{P}}_z - \frac{1}{4 \pi |z|^3} (I-3 e_z e^t_z) \int_{\{y:|y| \leq |z| \}} K(y, 2 \nu s) dy  \Big) 
g(x-z,t-s) \Bigg\} dz ds, \nonumber
\end{eqnarray}
then \(u\) is a weak solution to (\ref{ns}).
\end{proposition}

This integral formulation incorporating incompressibility
is key to the proof of Theorem \ref{thm1} and the other existence and uniqueness results given in Sections 4 and 5.  The following lemmas are used in the proof of the proposition.  Throughout the Fourier transform on \({\bf{R}}^3\) is defined as
\[
\hat{f} ( \xi ) = (2 \pi)^{-3/2} \int_{{\bf{R}}^3} e^{-ix \cdot \xi} f(x) dx
\]
and understood in a distributional sense.  Then also
\[
f(x) = (2 \pi )^{-3/2} \int_{{\bf{R}}^3} e^{ix \cdot \xi} \hat{f} ( \xi) d \xi.
\]

The bilinear forms \({\bf b}_1\) and \({\bf b}_2\) defined above and the matrix \(I-3e_ze^t_z\) can be bounded as follows.
\begin{lemma}\label{lembivop}  For any \(y \in {\bf R}^3\) with \(|y| \neq 0\) and \(u,v \in {\bf R}^3\),
 \begin{equation}\label{eqb1bd}
| {\bf b}_1(y;u,v) | \leq |u| \, |v|,
\end{equation}
\begin{equation}\label{eqb2bd}
|{\bf b}_2(y;u,v)| \leq 2 |u| \, |v|,
\end{equation}
and
\begin{equation}\label{eqcbd}
|(I- 3 e_y e_y^t ) u| \leq 2 |u|.
\end{equation}
\end{lemma}
\noindent{\it Proof: }  Fix \(y,u,v \in {\bf R}^3\) with \( |y| \neq 0 \).  Let \( \alpha = e_u \cdot e_y \), \(\beta = e_v \cdot e_y\), and \( \gamma = e_u \cdot e_v\).  For \(a,b \in [0,1]\), \( f( a , b ) := 
a(1 - b^2)^{1/2} + b(1- a^2)^{1/2} \leq 1\) with equality for \(a^2 + b^2=1\).
The triangle inequality gives 
\[
|{\bf b}_1(y;u,v)| =|u \, |v| \,|\alpha {\bf P}_y e_v + \beta {\bf P}_y e_u|
\leq |u| \, |v| f(| \alpha|, | \beta|) \leq |u | \, |v|.
\]
To derive (\ref{eqb2bd}), first notice that
\[
|{\bf b}_1(y;u,v)|^2=|u|^2|v|^2(\alpha^2(1-\beta^2)+2 \alpha \beta(\gamma-\alpha \beta)+\beta^2(1-\alpha^2))
\]
and
\[
| \gamma - \alpha \beta | = |e_u \cdot {\bf P}_y e_v| \leq |{\bf P}_y e_u| \, | {\bf P} e_v|
= (1 - \alpha^2)^{1/2} (1- \beta^2 )^{1/2} .
\]
Then
\begin{eqnarray*}
|{\bf b}_2(y; u,v)|^2 & = & |{\bf b}_1(y;u,v)|^2 +( u \cdot (I-3 e_y e_y^t )v)^2 \\
& = & |u|^2 \, |v|^2 \, ( \alpha^2 + \beta^2 + \alpha^2 \beta^2 + ( \gamma - 2 \alpha \beta)^2) \\
& \leq & |u|^2 \, |v|^2 \, ( \alpha^2 + \beta^2 + \alpha^2 \beta^2 + ( (1- \alpha^2)^{1/2}
(1- \beta^2)^{1/2} + | \alpha \beta|)^2) \\
& = & |u|^2 \, |v|^2 \, ( 1 + \alpha^2 \beta^2 +  2 |\alpha \beta | f( | \alpha|, (1- \beta^2)^{1/2})) \\
& \leq & |u|^2 \, |v|^2 \, ( 1 + | \alpha \beta | )^2 \\
& \leq & (2 |u| \, |v|)^2 .
\end{eqnarray*}
To verify (\ref{eqcbd}) note that
\[
|(I-3 e_y e_y^t )u|^2 = u \cdot (I + 3 e_y e_y^t) u = |u|^2 (1+3 \alpha^2) \leq 4|u|^2.
\]
\hfill \(\Box\)

\begin{lemma}\label{lemGamma}
Suppose that for each fixed \(s>0\), \(\Gamma ( \cdot ,s): {\bf{R}}^3 \to {\bf{R}}^{3 \times 3} \)
has Fourier transform
\[
\hat{\Gamma}(\xi,s) = (2 \pi)^{-3/2} e^{- \nu | \xi |^2 s} {\bf{P}}_\xi.
\] 
Then for \(x \in {\bf{R}}^3\) with \(|x| \neq 0\), 
\[
\Gamma (x,s) =K(x, 2 \nu s) {\bf{P}}_x -(4 \pi)^{-1} |x|^{-3} (I -3 e_x e^t_x ) \int_{\{y:|y| \leq |x|\}} 
K(y, 2 \nu s)dy.
\]
\end{lemma}
The result presented in Lemma \ref{lemGamma} is not new; it goes back at least to Solonnikov \cite{S64}, who used it to calculate estimates of solutions of the linearized Navier-Stokes equations.  It appears again in Koch and Solonnikov \cite{KS02} in a study of the Stokes problem.  More recently, 
Thomann and Guenther \cite{TG04} use it to compute an explicit formula for the fundamental solution of the linearized Navier-Stokes equation in terms of special functions.  The representation given here is more explicitly stated in terms of  \(K\), the transition density for Brownian motion, than that of Solonnikov \cite{S64}, Koch and Solonnikov \cite{KS02} or Thomann and Guenther \cite{TG04}.  The derivation is included for completeness.

\noindent{\it Proof: }  Fix \(s>0\) and notice that 
\(
\hat{K}( \xi,s) =(2 \pi)^{-3/2} e^{- \nu |\xi|^2s}.
\)
Suppose that \(\hat{\gamma}_0 ( \xi) = (2 \pi)^{-3/2}| \xi |^{-2} e^{-\nu | \xi |^2 s}.
\)
Some computation gives 
\begin{equation}\label{eqG1}
\gamma_0 (x) 
 =  (2 \pi)^{-3/2} (2 \nu s)^{-1/2} |x|^{-1} \int_{ \rho = 0}^{|x|} e^{- \rho^2 / 4 \nu s} d \rho.
\end{equation}
(For computational details of the Fourier transform for radial functions, see Folland \cite{F92}, page 247.)
If \(\gamma_{i,j}:{\bf{R}}^3 \to {\bf{R}}\) is defined via
\[
\hat{\gamma}_{i,j} ( \xi) = -(2 \pi)^{-3/2} \frac{ \xi_i \xi_j}{| \xi|^2} e^{- \nu | \xi |^2 s},
\]
then
\begin{eqnarray}\label{eqG2}
\gamma_{i,j}(x) & = & \frac{\partial^2 }{\partial x_i \partial x_j}\gamma_0 (x)  \nonumber \\
& = & (2 \pi )^{-3/2} (2 \nu s)^{-1/2} \bigg\{ \delta_{i,j} \Big(|x|^{-2} e^{-|x|^2/4 \nu s} 
-|x|^{-3} \int_{\rho =0}^{|x|} e^{- \rho^2 /4 \nu s} d \rho \Big)  \nonumber \\
& &  \quad +  x_i x_j |x|^{-2} \Big( 3 |x|^{-3} \int_{\rho = 0}^{|x|} e^{-  \rho^2 /4 \nu s} d \rho 
-(3|x|^{-2} + (2 \nu s)^{-1}) e^{-|x|^2/4 \nu s} \Big)  \bigg\}  .
\end{eqnarray}
Combining (\ref{eqG1}) and (\ref{eqG2}) gives
\begin{eqnarray*}
\Gamma_{i,j}(x,s) & = & K(x,2 \nu s) \delta_{i,j}  + \gamma_{i,j} (x) \\
& = & K(x, 2 \nu s) ({\bf{P}}_x)_{i,j}  \\
& & - (2 \pi )^{-3/2} (2 \nu s)^{-1/2} |x|^{-3} \;
\Big(\int_{\rho =0}^{|x|} e^{- \rho^2 / 4 \nu s} d \rho - |x| e^{-|x|^2 /4 \nu s}\Big) \, 
\Big( \delta_{i,j} - 3 x_i x_j |x|^{-2}\Big)
\end{eqnarray*}
To complete the proof note that 
\begin{equation}
\int_{\{y:|y| \leq |x| \}} K(y, 2 \nu s) dy 
 =  ( \pi \nu s )^{-1/2} \Big( \int_{\rho =0}^{|x|} e^{- \rho^2 / 4 \nu s} d \rho -|x|e^{-|x|^2/4 \nu s} \Big).
\end{equation}
\hfill \(\Box\)

The following lemma is given without proof. 
\begin{lemma}\label{lemdiff}  Suppose that \(V:{\bf{R}}^3 \to {\bf{R}}^3\) has
Fourier transform 
\[\hat{V}(\xi)= i(2 \pi)^{-3/2} \int_{{\bf{R}}^3} \hat{v} ( \eta ) \, \xi \cdot \hat{v} (\xi - \eta ) d \eta .
\]
Then \(V\) is given by 
\[
V_k(x) = \sum_{j=1}^3 \frac{ \partial}{\partial x_j} (v_k(x) v_j(x)) .
\]
\end{lemma}
\begin{lemma}\label{lemconvol}  For \(V\) as defined in Lemma \ref{lemdiff},
\begin{eqnarray*}
\int_{{\bf{R}}^3}  & & \Gamma (y,s) V(x-y) dy   =  - \int_{{\bf{R}}^3} 
\Big\{ (4 \nu s)^{-1} |y| K(y, 2 \nu s) \, {\bf b}_1(y;v(x-y),v(x-y))  \\
& & +  \Big(|y|^{-1} K(y, 2 \nu s) - 3 (4 \pi )^{-1} |y|^{-4} \int_{\{z: |z| \leq |y| \}} K(z, 2 \nu s) dz \Big) \,
{\bf b}_2(y;v(x-y),v(x-y)) \Big\}dy .
\end{eqnarray*}
\end{lemma}
\noindent{\it Proof: }  Lemma \ref{lemGamma} gives
\[
\Gamma_{i,j}(y,s) = \delta_{i,j} g_1(|y|) +y_i y_j g_2(|y|)
\]
where for \( r > 0\)
\[
g_1(r) = (4 \pi \nu s)^{-3/2} e^{-r^2/4 \nu s} -(4 \pi)^{-1} r^{-3} \int_{\{y: |y| \leq r \}} K(y, 2 \nu s) dy
\]
and
\[
g_2(r) = 3(4 \pi)^{-1} r^{-5} \int_{\{y: |y| \leq r \}} K(y, 2 \nu s) dy -(4 \pi \nu s)^{-3/2} r^{-2} e^{-r^2/4 \nu s}.
\]
A term by term integration by parts leads to the vector-valued function
\begin{eqnarray*}
 \int_{{\bf{R}}^3} \Gamma (y,s) V(x-y) dy  & = &  \int_{{\bf{R}}^3} \Big\{ (2 \nu s)^{-1}K(y, 2 \nu s)
 \Big((e_y \cdot v(x-y))^2 \, y -y \cdot v(x-y) \, v(x-y)\Big) \\
 & + & g_2(|y|) \Big( 2 (y \cdot v(x-y))v(x-y) +(|v(x-y)|^2 - 5(e_y \cdot v(x-y))^2)y \Big) \Big\} dy \\
  & = & \int_{{\bf{R}}^3} \Big\{-(4\nu s)^{-1} |y| K(y, 2 \nu s) \, {\bf b}_1(y;v(x-y),v(x-y)) \\
 & & \quad \quad
 + |y| g_2(|y|) \, {\bf b}_2(y;v(x-y),v(x-y)) \Big\}dy.
 \end{eqnarray*}
\hfill \(\Box\)
 
 \noindent{\it Proof of Proposition \ref{propins}:}  The incorporation of Leray's projection into the Fourier transformed Navier-Stokes equations as seen in (\ref{eqFNS}) gives a convenient formulation
 for weak solutions \(u\) of  (\ref{ns}).  For fixed \(s\), apply the inverse Fourier transform to the Fourier Navier-Stokes equations (\ref{eqFNS}), letting 
 \(V\) represent the vector-valued function with \(V_k(x,s) = \sum_{j=1}^3 \frac{\partial}{\partial x_j} (u_k(x,s) u_j(x,s))\), and use Lemmas \ref{lemGamma}, \ref{lemdiff} and \ref{lemconvol} to see that 
 \begin{eqnarray*}
 u(x,t) & = &  \int_{{\bf{R}}^3} u_0(x-y)K(y, 2 \nu t) dy \\
 & & \quad + \int_{s=0}^t \Bigg\{ -\int_{{\bf{R}}^3}  \Gamma(y,s) V(x-y,t-s) dy 
 + \int_{{\bf{R}}^3}  \Gamma(y,s) g(x-y,t-s) dy \Bigg\} ds \\
 & = &  \int_{{\bf{R}}^3} u_0(x-y)K(y, 2 \nu t) dy \\
 & &  + \int_{s=0}^t \Bigg\{ \int_{{\bf{R}}^3} \bigg((4 \nu s)^{-1} |y|K(y,2 \nu s) \, 
 {\bf b}_1(y; u(x-y,t-s),u(x-y,t-s))\\
 & &  + \Big(|y|^{-1}K(y,2 \nu s) -3(4 \pi)^{-1} |y|^{-4} \int_{\{z:|z| \leq |y| \}}K(z,2 \nu s )dz \Big)\\
 & & \quad \quad \quad
 {\bf b}_2 (y; u(x-y,t-s),u(x-y,t-s)) \bigg)dy \\
 & & +  \int_{{\bf{R}}^3}  \Gamma(y,s) g(x-y, t-s)dy \Bigg\} ds.
 \end{eqnarray*}
 \hfill \(\Box\)

\section{Navier-Stokes majorizing kernel pairs.}

The majorizing kernel pairs defined below are used to control the magnitude of the forcing and the initial velocity in a way that allows the control of the magnitude of the velocity over time.  The first component of the pair, \(h\), dominates the velocity.  The second component, \(\tilde{h}\), dominates the forcing in a way tailored to maintain the velocity domination.  The conditions specificing the majorizing kernels allow construction of transition probabilities based on the integral representation of the Navier-Stokes equations given in Proposition \ref{propins} above.  These kernels play a role in physical space  analogous to the role played in Fourier space by the scaling functions of Le Jan and Sznitman \cite{LJS97} and the Fourier multiplier majorizing kernels introduced by Bhattacharya et al in \cite{BCDGOOTW03}.
\begin{definition}  The pair \((h, \tilde{h}): {\bf{R}}^3 \times {\bf{R}}^3 \to (0, \infty ] \times [0, \infty ]\)
is a majorizing kernel pair with constant pair \( (\gamma, \tilde{\gamma}) \in (0, \infty ) \times [0, \infty )\) if
\(h\) is lower semi-continuous and uniformly locally square integrable,
\begin{equation}\label{eqmk}
 \sup_{x \in {\bf{R}}^3 } \frac{ \int_{{\bf{R}}^3} h^2(x-y) |y|^{-2} dy}{h(x)} = \gamma < \infty ,
\end{equation}
and
\[
  \sup_{x \in {\bf{R}}^3 } \frac{ \int_{{\bf{R}}^3} \tilde{h}(x-y) |y|^{-1} dy }{h(x)} =  \tilde{\gamma} < \infty .
\]
If \( \gamma =\tilde{\gamma} =1 \), we say that \((h,\tilde{h})\) is a standard majorizing kernel pair.
\end{definition}
Notice that both \(h\) and \(\tilde{h}\) are finite a.e. with respect to Lebesgue 
measure and \(\tilde{h}(x) =0\) a.e. if and only if \(\tilde{\gamma}=0\).
A simple scaling argument shows that if a majorizing kernel pair
 \((h,\tilde{h})\) has constant pair \( ( \gamma, \tilde{\gamma})\) with  
\(\tilde{\gamma} > 0\), 
then \( ( \gamma^{-1}h, ( \gamma \tilde{\gamma})^{-1} \tilde{h})\) is a standard majorizing kernel pair.

Recall from potential theory that an {\it excessive} function \(h: {\bf R}^3 \to (0, \infty ] \) satisfies both
\[
  \sup_{x \in {\bf{R}}^3 ,t >0} \frac{ \int_{{\bf{R}}^3} h(y) K(x-y,2 \nu t) dy}{h(x)} \leq 1
\quad \mbox{and} \quad \lim_{t \downarrow 0}  \int_{{\bf{R}}^3} h(y) K(x-y,2 \nu t) dy = h(x) .
\]
See Bass \cite{B95} or Doob \cite{D84} for background.  Majorizing kernel pairs of particular interest are those with the first member \(h\) being excessive.

\begin{definition}  A majorizing kernel pair \((h, \tilde{h})\) is an excessive majorizing kernel pair if \(h\) is excessive.
\end{definition}
The initial data \(u_0\) and the forcing \(g\) will be compared, respectively to \( h\) and \(\tilde{h}\).
\begin{definition}  The pair \((u_0,g) = (u_0, \{g( \cdot, t) : t \geq 0 \})\) is \( (h, \tilde{h})\)-admissible if
\[
\sup_{x \in {\bf R}^3} \frac{|u_0(x)| }{h(x)} < \infty \quad \mbox{and} \quad
\sup_{x \in {\bf R}^3, t \geq 0 } \frac{|g(x,t)|}{\tilde{h}(x)} < \infty.
\]
\end{definition}
\noindent
{\bf Example 1:}  Let \(h_0 (x) = |x|^{-1}\) and let \(\tilde{h}(x) = \tilde{h}_0(|x|)\) with 
\( \int_{{\bf R}^3} \tilde{h}(x) dx = 4 \pi \int_{r=0}^\infty r^2 \tilde{h}_0 (r) dr < \infty\).  Then \( (h_0, \tilde{h})\) is an excessive
 majorizing kernel pair with constant pair \( ( \pi^3, 4 \pi^2 \int_{r=0}^\infty r^2 \tilde{h}_0 (r) dr ) \).  This can 
be seen by observing that \(h_0\) is uniformly locally square integrable and for \(x \neq 0\),
\[
\int_{{\bf{R}}^3 }|y|^{-1} K(x-y,t) dy =  |x|^{-1} \; P(|Z| < |x|/\sqrt{t}), 
\]
where \(Z\) is a standard normal r.v.,
\[
\int_{{\bf{R}}^3 }|y|^{-2}|x-y|^{-2} dy = \pi^3 |x|^{-1},
\]
and
\[
|x| \int_{{\bf{R}}^3 }|y|^{-1}\tilde{h}(x-y)dy = 4 \pi \int_{r=0}^\infty \min\{r,|x| \} \; r  \; \tilde{h}_0 (r) dr .
\]
The value on the left above increases to \( 4 \pi \int_{r=0}^\infty r^2 \tilde{h}_0 (r) dr = \int_{{\bf R}^3} \tilde{h}_0 (|y|)dy \) as 
\(|x| \to  \infty \).  As an example, let \(\tilde{h}_0(r) =r^{-2}(1+r)^{-(1+ \varepsilon)}\) for some fixed \(\varepsilon > 0\).  Then \((u_0,g)\) is \((h_0,\tilde{h})\)-admissible if \(\sup_{x \in {\bf R}^3} |x||u_0(x)| < \infty\) and 
\(\sup_{ x \in {\bf R}^3, t > 0} |x|^2(1+|x|)^{1 + \varepsilon} |g(x,t)| < \infty \).

The Fourier transform of the function \(h_0\) given here is the Le Jan-Sznitman scaling function
\(|\xi|^{-2}\) that arises naturally in the study of the Fourier Navier-Stokes equations.  See Le Jan and Sznitman \cite{LJS97} for the original study.  See also Bhattacharya et al \cite{BCDGOOTW03}  for refinement  of the Fourier space majorizing kernel method.  See also the remark following Example 4.

The following propositions illustrate some of the structure of the class of majorizing kernel pairs.
\begin{proposition}  Let \((h,\tilde{h})\) be a standard (excessive) majorizing kernel pair.

\parskip = 0.05 in
\noindent
(a)  Then
\[
(h( \cdot - \mu),\tilde{h}( \cdot - \mu)), \; \; \mu \in {\bf{R}},
\] 
\[
(\sigma h( \sigma \cdot ), \sigma^3 \tilde{h}( \sigma \cdot )), \;  \; \sigma >0 ,
\]
and
\[
(h(A \cdot), \tilde{h}(A \cdot)), \; \; \mbox{ A a 3 by 3 matrix with} \; \; A^tA=I,
\]
are also standard (excessive) majorizing kernel pairs.

\noindent
(b)  Let \(F\) denote a probability distribution function on \({\bf{R}}^3\).  Then
\[
(\int_{{\bf R}^3} h( \cdot -y)dF(y), \int_{{\bf R}^3} \tilde{h}( \cdot -y)dF(y))
\]
is also a (excessive) majorizing kernel pair with constant pair \(( \gamma, \tilde{\gamma}) \in (0,1] \times (0,1]\).
\end{proposition}
\noindent{\it Proof: }  Part (a) is easily checked  via the appropriate change of variables.  Part (b)
follows from an application of the Cauchy-Schwartz inequality followed by Fubini's Theorem. \hfill \(\Box\)

\parskip = 0.2 in
\begin{proposition}\label{propmajk2}   Let \(\{(h_j,\tilde{h}_j) : j \geq 1\}\) be a sequence of (excessive) majorizing kernel pairs with corresponding constant pairs \(\{(\gamma_j,\tilde{\gamma}_j) : j \geq 1\}\).

\parskip = 0.05 in
\noindent
(a)  Then \((h_1 \wedge h_2, \tilde{h}_1 \wedge \tilde{h}_2)\) is a (excessive) majorizing kernel pair
with constant pair \((\gamma, \tilde{\gamma}) \in (0, \gamma_1 \wedge \gamma_2] \times
[0,\tilde{\gamma}_1 \wedge \tilde{\gamma}_2 ]\).

\noindent
(b)  For any \(p \in (0,1)\)
\[
(h_1^p \, h_2^{1-p}, \tilde{h}_1^p \,  \tilde{h}_2^{1-p})
\] 
is a (excessive) majorizing kernel pair with constant pair  \( (\gamma,\tilde{\gamma})  \in (0, \gamma_1^p \gamma_2^{1-p}] \times
[0, \tilde{\gamma}_1^p \tilde{\gamma}_2^{1-p} ]\).

\noindent
(c)  For \(\{p_j: j \geq 1\}\) with \(\sum_{j=1}^\infty p_j =1\) and \(p_j \geq 0\) for each \(j \geq 1\)
\[
(\sum_{j=1}^\infty p_j h_j , \sum_{j=1}^\infty p_j \tilde{h}_j )
\]
is a (excessive) majorizing kernel pair with constant pair \(  (\gamma,\tilde{\gamma}) \in 
(0, \sum_{j=1}^\infty p_j \gamma_j ] \times
[0, \sum_{j=1}^\infty p_j \tilde{\gamma}_j]\).
\end{proposition}
\noindent{\it Proof: }  Part (a) is obvious.  Apply H\"{o}lder's inequality to deduce (b).  Part (c) again follows from the Cauchy-Schwartz inequality and Fubini's Theorem. \hfill \(\Box\)

\parskip = 0.2 in
\begin{definition}
The functions \(h_j: {\bf R}^3 \to (0, \infty]\), \(j=1,2\) are equivalent if for some \( c \in (1,\infty)\),
\[
c^{-1} h_1(x) \leq h_2(x) \leq c h_1(x) \; \; for \; \; all \; \; x \in {\bf R}^3.
\]
\end{definition}
\begin{proposition}\label{propequiv}
If \(h_1\) and \(h_2\) are equivalent and \((h_1, \tilde{h})\) is a majorizing kernel pair, then \((h_2,\tilde{h})\) is also a majorizing kernel pair.
\end{proposition}
\noindent{\it Proof: }  The proof follows easily from the definitions.  \hfill \(\Box\)

\noindent
{\bf Example 2:}  The majorizing kernel \(h_0\) of Example 1 has a singularity at the origin.   Propositions 3.1 and 3.2, can be used to construct a kernel with a countable number of singularities as follows.  Take \(\{ \mu_j:j \geq 1 \}\) to be a sequence in \({\bf R}^3\) and \(\{p_j :j \geq 1\}\) with \(\sum_{j=1}^\infty p_j =1\) and \(p_j \geq 0\) for each \(j \geq 1\).  Let 
\[
h(x) = \sum_{j=1}^\infty p_j |x- \mu_j|^{-1} 
\]
and
\[
\tilde{h}(x) = \sum_{j=1}^\infty p_j \tilde{h}_0(|x-\mu_j|)
\]
for some \(\tilde{h}_0:[0,\infty) \to [0,\infty]\) with \(\int_{r=0}^\infty r^2 \tilde{h}_0(r)dr < \infty\).
Note that if \(\tilde{h}_0\) has a singularity, then \(\tilde{h}\) itself will have a countable number of singularities as well.  \((h,\tilde{h})\) is a majorizing kernel pair with constant pair \((\gamma, \tilde{\gamma}) \in (0,\pi^3] \times (0,  \int_{{\bf R}^3} \tilde{h}_0(|x|)dx]\).

\noindent
{\bf Example 3:}  The operations of Propositions 3.1 and 3.2 can also be used to produce bounded kernels based on \(h_0\).  For example, convolving \(h_0\) with the probability density 
\( f(y) = (2 \pi |y|)^{-1}(1+|y|)^{-3}\)
yields the bounded kernel \[
h_1 (x) = (1+|x|)^{-1}.
\]
Taking \(\tilde{h}(x) = \int_{{\bf R}^3} \tilde{h}_0(|x-y|) f(y) dy\) yields a majorizing kernel pair with constant pair \((\gamma, \tilde{\gamma}) \in (0,\pi^3] \times (0, \int_{{\bf R}^3} \tilde{h}_0(|x|)dx]\).

Using Proposition \ref{propequiv}  we see that any radial function \(h\) that decreases as \(|x|\) increases with \(\lim_{|x| \to 0} |x|h(x) =c \in [0,\infty)\) and \(\lim_{|x| \to \infty} |x| h(x) =C \in (0, \infty)\) can be used as a majorizing kernel when
paired with an appropriate \(\tilde{h}\).  Indeed kernels \(h\) of this description are in the Marcinkiewicz space \(L^{3,\infty}\); see Cannone and Karch \cite{CK03} for a related study of stability results for Navier-Stokes equations.

\noindent
{\bf Example 4:}  Let \(H(x) =(1 +|x|^2)^{-1}\).   Although \(H\) is not excessive, it is uniformly locally square integrable and 
\[
\int_{{\bf{R}}^3 } H(y) K(x-y,t) dy \leq (1+3e^{-2/3})H(x) \quad \mbox{for all} \quad  x \in {\bf R}^3.
\]
Additionally, 
\[
\int_{{\bf{R}}^3 } H^2(y) |x-y|^{-2} dy = \pi^2 H (x) \quad \mbox{for all} \quad x \in {\bf R}^3.
\]
However, it is not difficult to show that for any \(\tilde{h}\) with \(\int_{{\bf R}^3} \tilde{h}(y)dy >0\), 
\[
\sup_{x \in {\bf R}^3} (H(x))^{-1} \int_{{\bf{R}}^3 }|y|^{-1}\tilde{h}(x-y)dy = \infty.
\]
This gives the majorizing kernel pair \((H,0)\) with constant pair \((\pi^2,0)\).  Thus \((u_0,g)\) being \((H,0)\)-admissible
corresponds to \(\sup_{x \in {\bf R}^3} (1+|x|^2)|u_0(x)| < \infty\) and the forcing \(g\) being identically 0.
Note that \(\hat{H}(\xi) = \sqrt{\pi/2}|\xi|^{-1}e^{-|\xi|}\).  This is the second of the two Fourier Navier-Stokes majorizing kernels introduced by Le Jan and Sznitman \cite{LJS97}.

\noindent
{\bf Remark:}  Although Examples 1 and 4 are suggestive, it is not generally true that the Fourier transform of a majorizing kernel \(h\) satisfying (\ref{eqmk}) is a Fourier multiplier majorizing kernel or vice versa.  The key to the correspondence in Examples 1 and 4 is the equality
\[
\int h^2(x-y)|y|^{-2}dy = \gamma h(x).
\]
On the Fourier side, this becomes
\[
\hat{h}*\hat{h}(\xi)=c \gamma | \xi | \hat{h}(\xi)
\]
for a fixed constant \(c\).  The condition \(h*h(\xi) \leq C | \xi | h( \xi )\) is that required by Bhattacharya et al \cite{BCDGOOTW03} for Fourier multiplier majorizing kernels for the Fourier transformed Navier-Stokes equations.

\noindent
{\bf Example 5:}   For \(p \in (1,2]\), let \(H_p(x) = (1+|x|)^{-p}\).  Then \((H_p,0)\) is a majorizing kernel pair with constant pair \((\gamma,0)\), \(\gamma \leq \pi^{4-p}(1+1/\sqrt{2})^{p-1}\).  This can be deduced by first noting that \(H\) as defined in Example 4 is equivalent to \((1+|x|)^{-2}\) with \(c=1+1/\sqrt{2}\) and then applying (b) of Proposition \ref{propmajk2} using \(H\) and the kernel \((1+|x|)^{-1}\) of Example 3. 

\section{Construction of the stochastic representation and existence and uniqueness with excessive kernels.} 

This section opens by introducing the notation necessary to define tree-indexed Markov processes.  After defining some crucial transition densities, the integral equation (\ref{eqexpins}) of Proposition \ref{propins} is reformulated in terms of the transition densities and the velocity and forcing scaled respectively by \(h\) and \(\tilde{h}\) for a majorizing kernel pair \((h,\tilde{h})\).  The underlying tree-indexed Markov process and stochastic recursion are then defined and described.  The existence and uniqueness results for excessive kernels follow.

Let \({\cal{V}} := \cup_{n=0}^\infty \{0,1 \}^n\) denote the full binary tree with root \(\phi = \{0,1 \}^0\).  Let 
\(\partial {\cal{V} }:=  \{0,1 \}^{\bf{N}}\) denote the boundary of \(\cal{V}\).  For \(v=<v_1 ,..., v_n> 
\in \{ 0,1 \}^n\), we say that the magnitude of \(v\), \(|v|=n\).  For \(v \in \cal{V}\) 
define \(v|0= \phi \) and for \(v \neq \phi\) and \(n \geq 1\), 
\(v|n = <v_1,...,v_n>\).  If for some \(n < |v|\), \(v|n = w\), we say that \(w\) is an ancestor of \(v\) and, conversely, that \(v\) is a descendant of \(w\).  For \(v \in \cal{V}\) with \(|v|=n\), \(0 \leq n < \infty\), set \(\bar{v}= v|(|v|-1)\) and \(v \ast k =<v_1,...,v_n,k>\), \(k=0,1\).  Thus \(v\) is the child of \(\bar{v}\) and has children \(v \ast 0\) and \(v \ast 1\).  Although it may currently seem nonsensical, it is convenient in the sequel to interpret \(\bar{\phi}\) as the precursor of the root \(\phi\).  
For \(v,w \in \cal{V}\) with \(v \neq w\), denote the last common ancestor of \(v\) and \(w\) as 
\(v \wedge w\).  This is defined as follows.  Let \(n_{v,w} =  \max \{ m \geq 0 : v|m = w|m\}\) and then take \(v \wedge w = v|n_{v,w}\).
We say that \({\cal{W}} \subset \cal{V}\) is a {\it rooted binary sub-tree} if \(\phi \in \cal{W}\), for any \(v \in { \cal W}\), \(v|k \in {\cal W}\) for all \(k < |v|\), and for any \(v \ast j \in \cal{W}\), \(j=0,1\), then \( v \ast (1-j) \in \cal {W}\) as well.  For a finite rooted binary sub-tree \(\cal W\), define the boundary of \(\cal{W}\), \(\partial \cal{W} \), as the elements of \(\cal{W}\) which have no descendants in \(\cal{W}\); that is
\[
\partial { \cal{W}} = \{ v \in {\cal{W}} : v \ast 0 \not\in \cal{W} \}.
\]
The interior of \(\cal{W}\), \({\cal{W}}^\circ = {\cal{W}} \setminus \partial \cal{W} \), consists of the elements of \(\cal{W}\) that have descendants in \(\cal{W}\).  For example, a Galton-Watson tree with a single progenitor and offspring distribution concentrated on 0 and 2 is a rooted binary sub-tree.  If the expected number of offspring is less than or equal to 1, then the binary sub-tree is finite with probability 1.

Let \(( {\bf{ B}}, \cal{B})\) be a measurable space and let \({{\bf X}} = \{ X_v : v \in \cal{V} \} \) be a \(\cal{V}\)-indexed collection of \(\bf B\)-valued random variables defined on a common probability space.  Let \({\cal{F}}_v = \sigma(X_v) \), \({\cal{G}}_v = \sigma ( \{ X_{v|n} : n \leq |v| \} ) \) denote the \(\sigma\)-field generated by the collection of r.v.'s indexed by \(v\) and \(v\)'s ancestors, and 
\({\cal{H}}_v = \sigma ( \{ X_{w} : w \Big| |v| =v \} ) \) denote the \(\sigma\)-field generated by the collection of r.v.'s indexed by \(v\) and \(v\)'s descendants.  Let \({\cal F}_{\bar{\phi}}\) denote the trivial \(\sigma\)-field.
\begin{definition}  \(\bf{X}\) is a \(\cal{V}\)-indexed Markov process if for any \(v \in \cal{V}\) with \(|v| < \infty\),  \({\cal{H}}_{v \ast 0}\) and \({\cal{H}}_{v \ast 1}\) are conditionally independent given \({\cal F}_v\), so that 
\[
P( A_1 \cap A_2 | {\cal F}_v) = P(A_1| {\cal F}_v) P(A_2 | {\cal F}_v) \quad a.s. \; \;  P  \quad \mbox{for} \quad
A_k \in {\cal{H}}_{v \ast k}, k=0,1,
\]
and for any \(v,w \in \cal{V}\)
and any \({\cal{H}}_w \)-measurable random variable \(Y\) with \(E|Y| < \infty\),
\[
E(Y \big| {\cal{G}}_v) = E(Y \big| {\cal{F}}_{v \wedge w}) \quad a.s. \; P.
\]
\end{definition}
The distribution of a \(\cal{V}\)-indexed Markov process is completely specified by the conditional distributions of the \(X_{v}\) given \({\cal{F}}_{\bar{v}}\) for \(v \in \cal{V}\).

The integral form of the incompressible Navier-Stokes equation will be expressed in terms of an expectation of a functional of a \(\cal{V}\)-indexed Markov process.  The following conditional densities involving the majorizing kernel pairs will be used to specify the transition probabilities and thus the distribution.  
Let \((h,\tilde{h})\) be a majorizing kernel pair.  Fix \(x \in{\bf{R}}^3\) and define a probability density on
\( {\bf{R}}^3 \times {\bf{R}}^3\) via
\begin{equation}\label{equspdens}
f(y,z|x) =  \frac{|y|^{-1} |z|^{-4} h^2(x-z){\bf{1}}[|z|>|y|]}{2 \pi \int_{{\bf{R}}^3} |z|^{-2} h^2(x-z)dz}.
\end{equation}
If \(\tilde{h}\) is not identically 0, likewise define 
\begin{equation}\label{eqfspdens}
\tilde{f}(y,z|x) =  \frac{|y|^{-1} |z|^{-3} \tilde{h}(x-z){\bf{1}}[|z|>|y|]}{2 \pi \int_{{\bf{R}}^3} |z|^{-1}
 \tilde{h}(x-z)dz}.
\end{equation}
The above densities will provide spatial transition densities for our \(\cal V\)-indexed Markov process.  
Define two conditional waiting time densities for \(s>0\) and \(y \in {\bf{R}}^3\) via
\begin{equation}\label{eqwt0dens}
f_0(s|y) = (2 \pi)^{-1/2}(2 \nu)^{-3/2}s^{-5/2}|y|^3 e^{-|y|^2/4 \nu s} = 2 \pi s^{-1}|y|^3 K(y,2 \nu s)
\end{equation}
and
\begin{equation}\label{eqwt1dens}
f_1(s|y) = (4 \pi \nu)^{-1/2}s^{-3/2}|y| e^{-|y|^2/4 \nu s} = 4 \pi \nu |y| K(y,2 \nu s).
\end{equation}
Adjoin a trap state \(\theta\) to \({\bf R}^3\), and define 
\begin{equation}\label{eqhbmdens}
J(y,t|x) = \left\{
 \begin{array}
{l l}
(h(x))^{-1}{h(y) K(x-y,2 \nu t)} \quad & \mbox{if} \quad x,y \in {\bf R}^3, \; h(x) < \infty; \\
1- (h(x))^{-1}{\int_{{\bf R}^3}h(y) K(x-y,2 \nu t)dy} \; & \mbox{if} \quad x \in {\bf R}^3, y = \theta; \\
1 \quad & \mbox{if} \quad x = y = \theta \\
0 & \mbox{otherwise.}
\end{array} \right.
\end{equation}
If \(h\) is excessive, this is the transition density of a \(h\)-Brownian motion on \({\bf R}^3 \cup \{ \theta \}\).

For a majorizing kernel pair \((h, \tilde{h})\) and \(x \in {\bf R}^3 \cup \{ \theta \} \), define
\begin{equation}\label{eqchi0def}
\chi_0(x) = \left\{
 \begin{array}
{l l}
(h(x))^{-1} u_0(x)  \quad & \mbox{if} \quad x \in {\bf R}^3 \\
0  \quad & \mbox{if} \quad x= \theta
\end{array} \right.
\end{equation}
and
\begin{equation}\label{eqphidef}
\varphi(x,t) =\frac{g(x,t)}{\tilde{h}(x)}.
\end{equation}
Also let
\begin{equation}\label{eqmdef}
m(x) = \frac{\int_{{\bf R}^3} |y|^{-2} h^2(x-y)dy}{8 \pi \nu h(x)} ,
\end{equation}
and
\begin{equation}\label{eqmtdef}
\tilde{m}(x) = \frac{\int_{{\bf R}^3} |y|^{-1} \tilde{ h}(x-y)dy}{8 \pi \nu h(x)} .
\end{equation}

\begin{proposition}\label{propscins}
Suppose that the majorizing kernel pair \((h, \tilde{h})\) is used to define \(f\), \(\tilde{f}\), \(f_0\), \(f_1\), \(J\), \(\chi_0\), \(\varphi\), \(m\), and \(\tilde{m}\) as given in (\ref{equspdens})
through (\ref{eqmtdef}).  If 
\( h(x) \chi (x,t)\) is locally square integrable and \(\chi\) satisfies
\begin{eqnarray}\label{eqscins}
\chi(x,t) & = & \int_{{\bf R}^3} \chi_0(y) J(y,t|x) dy \\
\nonumber & + & \int_{s=0}^t \int_{y \in {\bf R}^3}  \int_{z \in {\bf R}^3} \bigg\{m(x) 
\bigg(f_0(s|z)f(y,z|x) {\bf b}_1(z; \chi(x-z,t-s),\chi(x-z,t-s)) \\
\nonumber & & \quad +  \Big(2f_1(s|z)f(y,z|x) -3f_1(s|y)f(y,z|x)\Big){\bf b}_2(z; \chi(x-z,t-s),\chi(x-z,t-s)) \bigg) \\
\nonumber & & \quad \quad +  \tilde{m}(x)  \Big(2f_1(s|z) \tilde{f}(y,z|x) {\bf P}_z \varphi(x-z,t-s) \\
\nonumber & & \quad \quad \quad -  f_1(s|y) \tilde{f}(y,z|x) (I-3 e_z e_z^t) \varphi (x-z,t-s)\Big) \bigg\} dz dy ds
\end{eqnarray}
then \(u(x,t)= h(x) \chi(x,t)\) is a weak solution to (\ref{ns}).
\end{proposition}

\noindent{\it Proof: }  Suppose that \(\chi(x,t)\) satisfies (\ref{eqscins}) above and define \(u(x,t) = h(x) \chi(x,t)\).  Note that
\[
\int_{y \in{\bf R}^3} f(y,z|x)dy = \frac{|z|^{-2} h^2(x-z)}{\int_{{\bf R}^3}|y|^{-2}h^2(x-y)dy}
\]
and
\[
\int_{y \in{\bf R}^3} \tilde{f}(y,z|x)dy = \frac{|z|^{-1} \tilde{h}(x-z)}{\int_{{\bf R}^3}|y|^{-1} \tilde{h}(x-y)dy}.
\]
Then
\begin{eqnarray*}
u(x,t) & = & \int_{y \in{\bf R}^3}  u_0(y) K(x-y,2 \nu t) dy \\
&+&\int_{s=0}^t \int_{z \in {\bf R}^3} 
\Bigg\{ \frac{h(x)m(x)}{\int_{{\bf R}^3}|y|^{-2}h^2(x-y)dy}
\bigg( \frac{ f_0(s|z)}{|z|^2} \; {\bf b}_1(z;u(x-z,t-s),u(x-z,t-s)) \\
& & \quad \quad + \;  \Big(\frac{2 f_1(s|z)}{|z|^2} - \int_{|y| \leq |z|} \frac{3f_1(s|y)}{2 \pi |y| |z|^4} dy \Big) \;
{\bf b}_2(z;u(x-z,t-s),u(x-z,t-s)) \bigg)  \\
& & \quad +  \frac{h(x) \tilde{m} (x)}{\int_{{\bf R}^3}|y|^{-1} \tilde{h}(x-y)dy} 
\bigg( \frac{2 f_1(s|z)}{|z|} \; {\bf P}_z g(x-z,t-s) \\
& & \quad \quad \quad \quad -  \int_{|y| \leq |z|} \frac{f_1(s|y)}{2 \pi |y| |z|^3} dy \; (I-3 e_z e_z^t)g(x-z,t-s)\bigg) \Bigg\} dz ds \\
& = & \int_{y \in{\bf R}^3}  u_0(y) K(x-y,2 \nu t) dy \\
&+&\int_{s=0}^t \int_{z \in {\bf R}^3} 
\Bigg\{ \bigg( \frac{|z| K(z,2 \nu s) }{4 \nu s} \; {\bf b}_1(z;u(x-z,t-s),u(x-z,t-s)) \\
& & \quad \quad + \;  \Big(\frac{K(z,2 \nu s) }{|z|} - \frac{3}{4 \pi |z|^4} \int_{|y| \leq |z|} K(y,2 \nu s) dy \Big) \;
{\bf b}_2(z;u(x-z,t-s),u(x-z,t-s)) \bigg)  \\
& & \quad +  
\bigg( K(z,2 \nu s) \; {\bf P}_z g(x-z,t-s) \\
& & \quad \quad \quad \quad - \frac{1}{4 \pi |z|^3} \int_{|y| \leq |z|} K(y,2 \nu s) dy \; (I-3 e_z e_z^t)g(x-z,t-s)\bigg) \Bigg\} dz ds.
\end{eqnarray*} 
The result then follows from Proposition \ref{propins}.
 \hfill \(\Box\)

It is helpful to combine the spatial and temporal transition densities \(f\), \(\tilde{f}\), \(f_0\) and \(f_1\) given in equations (\ref{equspdens}) through (\ref{eqwt1dens}) via the following randomization.  Let \(\{p_k : 1 \leq k \leq 5 \}\) satisfy \(p_k>0\) and \(\sum_{k=1}^5 p_k =1\).  
To appropriately balance the randomization, fix \(p \in (0, 1/2]\) and take
 \(p_1 = p/11\), \(p_2 = 4p/11\), \(p_3 = 6p/11\),
and \(p_4 = p_5 = (1-p)/2\).  Note that \(p_1 +p_2+p_3 =p \leq 1/2\).  For fixed \(x \in {\bf R}^3\) let
\begin{equation}\label{eqjtdens}
f(y,z,s,k|x) = \left\{
 \begin{array}
{l l}
p_1 f_0(s|z)f(y,z|x) \quad & \mbox{if} \quad k=1 \\
p_2 f_1(s|z)f(y,z|x) \quad & \mbox{if} \quad k=2 \\
p_3 f_1(s|y)f(y,z|x) \quad & \mbox{if} \quad k=3 \\
p_4 f_1(s|z) \tilde{f}(y,z|x) \quad & \mbox{if} \quad k=4 \\
p_5 f_1(s|y)\tilde{f}(y,z|x) \quad & \mbox{if} \quad k=5.
\end{array} \right.
\end{equation}
This can be thought of as the joint conditional density of a quadruple \((Y, Z, \tau, \kappa)\) taking values in \({\bf R}^3 \times {\bf R}^3 \times (0, \infty ) \times \{1,2,3,4,5 \} \).

We now define the tree-indexed Markov process that underlies our probabilistic representation.  Let 
\({\bf X} = \{( \{V_v(t):t \geq 0 \}, X_v, Y_v,Z_v, \tau_v, \kappa_v) : v \in{\cal V}\}\) be a \(\cal V\)-indexed Markov process with the following transition probabilities.  First assume that, given that \(X_{\bar{v}} =x\), \(\{V_v (t): t \geq 0\}\) is a \(h\)-Brownian motion with initial value \(x\) and
transition density \(J\) as defined in (\ref{eqhbmdens}) above.  Note that each \(V_v\) is itself a Markov process which can be assumed to be continuous with probability 1; see Doob \cite{D84}, Part 2, Chapter X for background.  Take the transition density of \((Y_v,Z_v, \tau_v, \kappa_v)\) given that \(X_{\bar{v}} =x\) to be given by \(f(y,z,s,k|x)\) as defined in (\ref{eqjtdens}), and then set 
\[
X_v= X_{\bar{v}} -Z_v.
 \]
Finally assume that for all \(v \in \cal V\), the process \(\{V_{v}(t) : t \geq 0 \}\) and the ensemble 
\((X_{v}, Y_{v},Z_{v}, \tau_{v}, \kappa_{v})\) are conditionally independent given \(X_{\bar{v}}\). 

Notice that the distribution of \(( \{V_v(t) : t \geq 0 \}, X_v, Y_v,Z_v, \tau_v, \kappa_v)\) only depends on the ensemble 
\((\{V_{\bar{v}}(t) : t \geq 0 \}, X_{\bar{v}}, Y_{\bar{v}}, Z_{\bar{v}}, \tau_{\bar{v}}, \kappa_{\bar{v}}) \) through the r.v. \(X_{\bar{v}}\).  The ensemble \(\{X_v:v \in {\cal V}\}\) is itself a branching Markov random walk.  The distribution of \(\bf X\) then depends on the initial value \(X_{\bar{\phi}}\).   Denote the  probability measure corresponding to \(X_{\bar{\phi}}=x\) by \(P_x\) and the
expectation with respect to this probability measure by \(E_x\).

For \(v \in {\cal V}\), define the random functional \(\Upsilon_v\) as follows.  Let
\({\bf B}_v\) denote the random bilinear operator 
\begin{equation}\label{eqBdef}
{\bf B}_v ( \; \cdot \; , \; \cdot \; ) = \left\{ \begin{array}{l l}
{\bf b}_1(Z_v; \; \cdot \;  , \; \cdot \; )  \quad & \mbox{if} \quad \kappa_v =1 \\
(-1)^{\kappa_v}{\bf b}_2(Z_v; \; \cdot \;  , \; \cdot \; )/2 \quad & \mbox{if} \quad \kappa_v=2,3
\end{array} \right.
\end{equation}
and \({\bf C}_v\) denote the random matrix 
\begin{equation}\label{eqCdef}
{\bf C}_v = \left\{ \begin{array}{l l}
{\bf P}_{Z_v} \quad & \mbox{if} \quad \kappa_v = 4 \\
-(I-3 e_{Z_v} e_{Z_v}^t)/2 \quad & \mbox{if} \quad \kappa_v =5 .
\end{array} \right.
\end{equation}
Then set
\begin{eqnarray}\label{eqUdef2}
\Upsilon_v(t)   =   \chi_0(V_v(t)) 
& + & \frac{11m(X_{\bar{v}})}{p} \; {\bf B}_v(\Upsilon_{v \ast 0}(t-\tau_v), \Upsilon_{v \ast 1}(t-\tau_v))
{\bf 1}[\kappa_v =1,2,3] \cap [\tau_v \leq t]  \nonumber \\
& + & \frac{4 \tilde{m}(X_{\bar{v}})}{1-p} \; {\bf C}_v \varphi(X_v,t-\tau_v) {\bf 1}[\kappa_v =4,5] \cap [\tau_v \leq t]. 
\end{eqnarray}
Notice that the recursion implicit in the definition of \(\Upsilon_v\) terminates if \(\tau_v >t\) or \(\kappa_v =4\) or 5.  That is
\begin{equation}\label{eqU0}
\Upsilon_v(t) {\bf 1}[\tau_v > t] = \chi_0(V_v(t))
\end{equation}
and
\begin{equation}\label{eqU1}
\Upsilon_v(t) {\bf 1}[\kappa_v =4,5]{\bf 1}[\tau_v \leq t] =  \chi_0(V_v(t))
+ \frac{4 \tilde{m}(X_{\bar{v}})}{1-p} \; {\bf C}_v  \varphi(X_v,t-\tau_v).
\end{equation}

The random space-time branching mechanism giving rise to \(\Upsilon_\phi\) can be explained as follows.  Initialize the Markov process \(\bf{X}\) by taking \(X_{\bar{\phi}}=x  \in {\bf R}^3\).  Activate the \(h\)-Brownian motion \(V_\phi(t)\) starting at \(X_{\bar{\phi}}\) at time \(t=0\).  This motion will run until time \(t\) with \(\chi_0(V_\phi (t))\) being used in calculating the first term of \(\Upsilon_\phi(t)\).  If \(V_\phi\) has entered the trap state \(\theta\) this term has value 0.  Although the process \(V_\phi\) runs until time \(t\), if \(\tau_\phi\) is less than \(t\) the path is deactivated at time \(\tau_\phi\).  If deactivation occurs before time \(t\), there are two possibilities.  If \(\kappa_\phi = 4\) or 5, the forcing is input into
\(\Upsilon_\phi\) by evaluating \({\bf C}_\phi \varphi\) at location \(X_\phi\) and time \( t-\tau_\phi\).
If \(\kappa_\phi = 1, 2\) or 3, two new active paths \(V_0\) and \(V_1\) are started at location  \(X_\phi\) at time \(\tau_\phi\).  The distributions of these newly activated paths are conditionally independent given  \(X_\phi\).  
If the paths \(V_0\) and \(V_1\) have been activated, then the process repeats on each of them with \(X_\phi\) taking the place of  \(X_{\bar{\phi}}\) and \(t-\tau_\phi \) taking the place of \(t\).  This is illustrated in Figure 1.  Notice that all paths activated by time \(t\) are used in calculating  \(\Upsilon_\phi(t)\), even though they may be deactivated before time \(t\).  For example evaluating the stochastic recursion for Figure 1 gives
\[
\Upsilon_\phi(t)  =  \chi_0(V_\phi(t)) +\frac{11m(X_{\bar{\phi}})}{p} {\bf B}_\phi( \Upsilon_0(t-\tau_\phi), \Upsilon_1(t-\tau_\phi))
\] with
 \[ \Upsilon_0(t-\tau_\phi) = \chi_0(V_0(t-\tau_\phi)) + \frac{4 \tilde{m}(X_\phi)}{1-p} {\bf C}_0 \varphi(X_0,t- \tau_\phi - \tau_0 )
 \]
  and 
\[
\Upsilon_1(t-\tau_\phi) = \chi_0(V_1(t-\tau_\phi)) + \frac{11 m(X_\phi )}{p} {\bf B}_1( \chi_0(V_{1,0}(t-\tau_\phi -\tau_1 )) , \chi_0(V_{1,1}(t-\tau_\phi -\tau_1))).  
\]
Here \(t\) is evaluated at the end of the time interval illustrated.

\begin{figure}\label{fig1}
\begin{center}
\scalebox{0.90}{\includegraphics{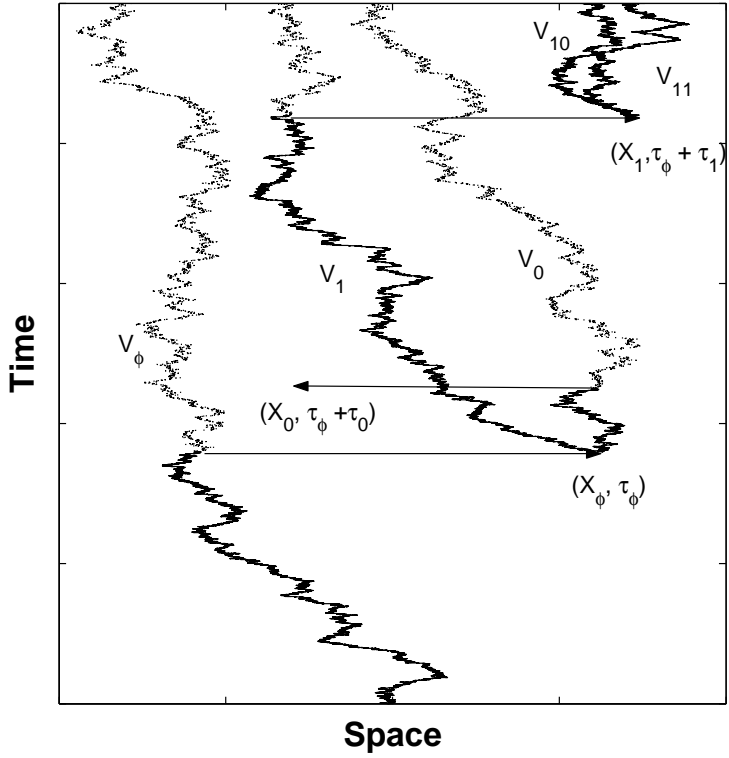}}
\end{center}
\caption{An illustration of an initial segment of the random process \(\bf X\) with \(\kappa_\phi =1,2\) or 3, \(\kappa_0 =4\) or 5, and \(\kappa_1 =1,2\) or 3.  The deactivated paths are represented by dotted lines.  }
\end{figure}

\noindent
{\bf Remark:}  Although it is not obvious from the construction above, the waiting times \(\tau_v\) can be tied to first passage times for Brownian motion.  Here's a sketch of an alternate construction that highlights the relationship.  For simplicity, just consider a single \(v \in \cal V\).  Let \(\{W(t): t \geq 0\}\) be a Brownian motion in \({\bf R}^3\) with \(W(0)=0\) and transition density \(K(y, 2 \nu t)\), independent of \(X_{\bar{v}}\) and the 
ensemble \((\{V_v(t): t \geq 0\},X_v,Y_v,Z_v,\kappa_v)\).  Write \(W(t)=(W^{(1)}(t),W^{(2)}(t),W^{(3)}(t))\)
and for \(a>0\) let
\[
T^{(i)}(a) = \inf \{ s>0: W^{(i)}(s)=a\} \quad \mbox{for} \quad i=1,2,3.
\]
The \(T^{(i)}(a)\)'s are independent stopping times 
with individual densities \((4 \pi \nu)^{-1/2}s^{-3/2}a e^{-a^2/4 \nu s}\).  
Let \[
T^{(0)}(a) = ( \sum_{i=1}^3 (T^{(i)}(a))^{-1})^{-1}.\]
It is easy to check that \(T^{(0)}(a)\) has density \((2 \pi)^{-1/2}(2 \nu)^{-3/2} s^{-5/2} a^3 e^{-a^2/4 \nu s}\). 
Note that \(T^{(0)}\) is not a stopping time with respect to the usual filtration associated with the Brownian motion \(W\).  Let
\[
 \tilde{\tau}_v= \left\{
 \begin{array}
{l l}
T^{(0)}(|Z_v|) \quad & \mbox{if} \quad \kappa_v = 1\\
T^{(1)}(|Z_v|) \quad & \mbox{if} \quad \kappa_v = 2,4\\
T^{(1)}(|Y_v|) \quad & \mbox{if} \quad \kappa_v = 3,5.\\
\end{array} \right.
\]
It is easy to check that, conditioned on the value of \(X_{\bar{v}}\), the joint distribution of 
\((\{V_v(t): t \geq 0\},X_v,Y_v,Z_v,\tilde{\tau}_v,\kappa_v)\),
 is the same as that of 
\((\{V_v(t): t \geq 0\},X_v,Y_v,Z_v, \tau_v,\kappa_v)\).  
This allows us to think of the \(\tau_v\)'s as waiting times that reflect excursion levels in orthogonal directions of a hidden independent Brownian motion in \({\bf R}^3\).

\begin{theorem}\label{thmbdsol}  Let \((h, \tilde{h})\) be an excessive majorizing kernel pair.  
If \(\Upsilon_\phi\) is as defined in (\ref{eqUdef2}) and 
\(
E_x| \Upsilon_\phi(t)| \leq M < \infty \) for all \(x\) and \(t\), then 
\[ u(x,t) = h(x) E_x \Upsilon_\phi(t)\]
is a weak solution to the Navier-Stokes equations.
\end{theorem}

\noindent{\it Proof: }  \(\Upsilon_0(t - \tau_\phi)\) and \(\Upsilon_1(t - \tau_\phi)\) are conditionally independent given \({\cal F}_\phi\), giving
\begin{eqnarray*}
 & & E_x( \Upsilon_\phi (t) | {\cal F}_\phi)  =   \chi_0(V_\phi (t)) \\
& & \quad \quad  +  m(x) \bigg\{ p_1^{-1}
{\bf b}_1(Z_\phi; E_x(\Upsilon_0( t - \tau_\phi)|{\cal F}_\phi ),
 E_x(\Upsilon_1( t - \tau_\phi)|{\cal F}_\phi)) {\bf 1}[ \kappa_\phi =1] \\
& & \quad \quad +  \Big( 2p_2^{-1} {\bf 1}[ \kappa_\phi =2] - 3( p_3)^{-1} 
{\bf 1}[\kappa_\phi = 3]  \Big) \\
& & \quad \quad \quad \quad
{\bf b}_2(Z_\phi; E_x(\Upsilon_0( t - \tau_\phi)|{\cal F}_\phi ),
 E_x(\Upsilon_1( t - \tau_\phi)|{\cal F}_\phi)) \bigg\} 
{\bf 1}[  \tau_\phi \leq t] \\
& & \quad \quad  +  \tilde{m} (x)\bigg\{ 2p_4^{-1} {\bf 1}[\kappa_\phi =4] {\bf P}_{Z_\phi}
 -(p_5)^{-1} {\bf 1}[ \kappa_\phi =5] (I - 3 e_{Z_\phi} e^t_{Z_\phi}) \bigg\}
  \varphi(x-Z_\phi,t- \tau_\phi){\bf 1}[ \tau_\phi \leq t] .
\end{eqnarray*}
Since also
\begin{equation}
E_x (\Upsilon_i(t-\tau_\phi)|{\cal F}_\phi){\bf 1}[Z_\phi =z,\tau_\phi =s] =E_{x-z}\Upsilon_\phi(t-s) 
\quad \hbox{for} \quad i=0,1,
\end{equation}
 \begin{eqnarray*}
E_x \Upsilon_\phi (t) & = & E_x ( E_x( \Upsilon_\phi (t) | {\cal F}_\phi) )\\
& = &  \int_{{\bf R}^3} J(y,t|x) \chi_0(y) dy \\
& + & \int_{s=0}^t  \int_{y \in {\bf R}^3}  \int_{z \in {\bf R}^3} \Bigg\{ m(x) \bigg(
f_0(s|z)f(y,z|x) {\bf b}_1(z; E_{x-z} \Upsilon_\phi(t-s),E_{x-z} \Upsilon_\phi(t-s)) \\
& + & \Big(2f_1(s|z)f(y,z|x) -3f_1(s|y)f(y,z|x) \Big)
{\bf b}_2(z; E_{x-z} \Upsilon_\phi(t-s),E_{x-z} \Upsilon_\phi(t-s)) \bigg)\\
& + & \tilde{m}(x)  \bigg(2f_1(s|z) \tilde{f}(y,z|x) {\bf P}_z \varphi(x-z,t-s) \\
& & \qquad - \; f_1(s|y) \tilde{f}(y,z|x) (I-3 e_z e_z^t) \varphi (x-z,t-s) \bigg) \Bigg\} dz dy  ds .
\end{eqnarray*}
Then \(u(x,t)\) is uniformly locally square integrable since \(|u(x,t)|=h(x)|E_x \Upsilon_\phi (t)| \leq Mh(x)\).  From Proposition \ref{propscins} the function
\(h(x) E_x \Upsilon_\phi (t)\) is a weak solution to (\ref{ns}).  \hfill \(\Box\)
\begin{theorem}\label{thmeu}  (Existence and Uniqueness)  Let \((h,\tilde{h})\) be an excessive majorizing kernel pair with constant pair \((\gamma, \tilde{\gamma}) \in (0, 8 \pi \nu p/11] \times [0, 2 \pi \nu (1-p)]\) for some \(p \in (0,1/2]\).  If \((u_0,g)\) is \((h,\tilde{h})\)-admissible with
\[
\sup_{x \in {\bf R}^3} \frac{|u_0(x)|}{h(x)} \leq \alpha \varepsilon \quad \mbox{and} \quad
\sup_{x \in {\bf R}^3, t \geq 0} \frac{|g(x,t)|}{\tilde{h}(x)} \leq \beta \varepsilon
\]
for some \(\alpha \in [0,1)\) and \(\varepsilon, \; \beta \in (0,1- \alpha)\),
then 
\(
u(x,t)= h(x) E_x \Upsilon_\phi (t)
\)
is a weak solution to the Navier-Stokes equations with
\[
\sup_{x \in {\bf R}^3, t \geq 0} \frac{|u(x,t)|}{h(x)} \leq \varepsilon.
\]
Furthermore, this solution is unique in the class 
\(
\{ v \in ({\cal S}' ({\bf R}^3 \times (0,\infty)))^3 :
 \sup_{x \in {\bf R}^3, t \geq 0} \frac{|v(x,t)|}{h(x)} \leq \varepsilon\}.
\)
\end{theorem}

The conditions of Theorem \ref{thmeu} are sometimes more easily checked in the following form.
\begin{corollary}\label{coreu}  Let \((h,\tilde{h})\) be an excessive majorizing kernel pair with constant pair \((\gamma, \tilde{\gamma})\), \(\tilde{\gamma} >0\),  and suppose that \((u_0,g)\) is \((h,\tilde{h})\)-admissible with
\[
\sup_{x \in {\bf R}^3} \frac{|u_0(x)|}{h(x)} \leq \frac{8 \pi \nu p \alpha \varepsilon}{11 \gamma}
 \quad \mbox{and} \quad
\sup_{x \in {\bf R}^3, t \geq 0} \frac{|g(x,t)|}{\tilde{h}(x)} \leq \frac{(4 \pi \nu )^2 p (1-p) \beta \varepsilon} { 11 \gamma \tilde{\gamma}} 
\]
for some \(p \in [0,1/2]\), \(\alpha \in [0,1)\) and \(\varepsilon, \; \beta \in (0,1- \alpha)\),
then 
\(
u(x,t)= h(x) E_x \Upsilon_\phi (t)
\)
is a weak solution to the Navier-Stokes equations with
\[
\sup_{x \in {\bf R}^3, t \geq 0} \frac{|u(x,t)|}{h(x)} \leq \frac{8 \pi \nu p \varepsilon }{11 \gamma }.
\]
This solution is unique in the class 
\(
\{ v \in ({\cal S}' ({\bf R}^3 \times (0,\infty)))^3 :
 \sup_{x \in {\bf R}^3, t \geq 0} \frac{|v(x,t)|}{h(x)} \leq \frac{8 \pi \nu p \varepsilon }{11 \gamma} \}.
\)
\end{corollary}
\noindent{\it Proof of Corollary \ref{coreu}:}  Rescale \((h,\tilde{h})\) and consider the majorizing kernel pair 
\[(\frac{8 \pi \nu}{11 \gamma} \, h, \frac{(4 \pi \nu )^2 p(1-p)}{11 \gamma \tilde{\gamma}} \,  \tilde{h})\] with constant pair \((8 \pi \nu p /11, 2 \pi \nu (1-p))\).  \hfill \(\Box\)

The next corollary treats the case of no forcing.
\begin{corollary}  Let \((h,0)\) be an excessive majorizing kernel pair with constant pair \((\gamma, 0)\) and suppose that 
\[
\sup_{x \in {\bf R}^3} \frac{|u_0(x)|}{h(x)} \leq \frac{4 \pi \nu \alpha \varepsilon}{11 \gamma}
\]
for some \(\alpha \in [0,1)\) and \(\varepsilon  \in (0,1- \alpha)\),
then 
\(
u(x,t)= h(x) E_x \Upsilon_\phi (t)
\)
is a weak solution to the Navier-Stokes equations with
\[
\sup_{x \in {\bf R}^3, t \geq 0} \frac{|u(x,t)|}{h(x)} \leq \frac{4 \pi \nu \varepsilon }{11 \gamma }.
\]
This solution is unique in the class 
\(
\{ v \in ({\cal S}' ({\bf R}^3 \times (0,\infty)))^3 :
 \sup_{x \in {\bf R}^3, t \geq 0} \frac{|v(x,t)|}{h(x)} \leq \frac{4 \pi \nu \varepsilon }{11 \gamma} \}.
\)
\end{corollary}
  
The proof of Theorem \ref{thmeu} depends upon Proposition \ref{propasbd} that follows. 
\begin{proposition}\label{propasbd}   Let \((h,\tilde{h})\) be an excessive majorizing kernel pair with constant pair
\((\gamma, \tilde{\gamma}) \in (0,8 \pi \nu \eta p/11] \times [0, 2 \pi \nu \eta (1-p)]\) for some \( p \in (0,1/2]\) 
and \(\eta > 0\).
 Suppose that for some \(\alpha \in [0,1)\) and \(\varepsilon, \beta \in (0, (1- \alpha)/\eta)\), 
\((u_0,g)\) is \((h,\tilde{h})\)-admissible with
\[
\sup_{x \in {\bf R}^3} \frac{|u_0(x)|}{h(x)} \leq \alpha \varepsilon
\quad
\mbox{and}
\quad
\sup_{x \in {\bf R}^3, t \geq 0} \frac{|g(x,t)|}{\tilde{h}(x)} \leq \beta \varepsilon .
\]
Then for all \(x \in {\bf R}^3\) and \(t>0\) 
\[
|\Upsilon_\phi(t)| \leq \varepsilon \quad a.s. \, \, P_x.
\]
\end{proposition}
The following lemma is used in the proof of the proposition.
\begin{lemma}\label{lemitcontract}  Let \({\cal W} \subset {\cal V}\) 
be a finite binary sub-tree.  Suppose that
\(\{{\bf b}_v:{\bf R}^3 \times{\bf R}^3 \to {\bf R}^3 :v \in {\cal W}\}\) has
\[
\sup_{v \in {\cal W}} |{\bf b}_v(x,y)| \leq |x| \, |y| \quad \mbox{for all} \quad x,y \in {\bf R}^3,
\]
and  \(\{(y_v,z_v,\eta_v,\sigma_v) \in {\bf R}^3 \times {\bf R}^3 \times [0, \infty) \times \{0,1\}: v \in {\cal W}\}\)
satisfies
\[
\sup_{v \in {\cal W}} \eta_v \leq \eta < \infty \]
and 
\[
\sup_{v \in  {\cal W}} |y_v| \leq \alpha \varepsilon\quad \mbox{and} \quad
\sup_{v \in \partial {\cal W}} |z_v| \leq \beta \varepsilon \quad \mbox{for some} \quad \alpha \in [0,1)
\quad \mbox{and} \quad  \varepsilon, \beta \in (0, (1- \alpha)/ \eta ].
\]
Then \(x_v\) defined iteratively on \(\cal W\) via
\[
x_v = \left\{ \begin{array}{l l}
y_v + \sigma_v \eta_v z_v \quad & \mbox{if} \quad v \in \partial {\cal W} \\
y_v + \eta_v {\bf b}_v(x_{v \ast 0},x_{v \ast 1}) \quad & \mbox{if} \quad v \in {\cal W}^\circ .
\end{array} \right.
\]
satisfies
\[
\sup_{v \in {\cal W}} |x_v| \leq \varepsilon .
\]
\end{lemma}
\noindent{\it Proof: }  
The result follows by induction starting on \(\partial {\cal W}\).  If \(v \in \partial \cal W\),
then
\[
|x_v|  \leq |y_v| +  \eta_v |z_v| \leq \alpha \varepsilon +\eta \beta \varepsilon \leq \varepsilon.
\]
If \( v \in {\cal W}^\circ\) and \(|x_{v \ast k}| \leq \varepsilon\) for \(k=0,1\), then
\[
|x_v| \leq |y_v| + \eta_v |{\bf b}_v(x_{v \ast 0},x_{v \ast 1})| \leq \alpha \varepsilon + \eta \varepsilon^2 
\leq \varepsilon.
\] \hfill \(\Box\)

\noindent{\it Proof of Proposition \ref{propasbd}:}  Let \({\cal W} = \{ v \in {\cal V} : \kappa_{v_j} =1,2,3 \quad \mbox{for all} \quad j <|v|\}\). 
Then \(\cal W \)  is a random binary sub-tree.  Indeed, \(\cal W \) corresponds to a Galton-Watson tree with each individual having either 0 or 2 offspring.  The probability of 2 offspring is \(p=P(\kappa_v =1,2,3) \leq 1/2\)
and the probability of 0 offspring is \(1-p=P( \kappa_v =4,5)\).   With probability 1, \(\cal W\) is finite; see for example Harris \cite{H63}.   
Set  \(X_{\bar{\phi}} = x \in  {\bf R}^3\) and \(S_\phi =0\).  For \(v \in \cal W\) with \(|v| > 0\), let \(S_v = \sum_{k=0}^{|v|-1} \tau_{v|k}\).  For
fixed \(t < \infty\), \(\Upsilon_\phi (t)\) is a functional of the random ensembles indexed by the nodes of the a.s. finite binary 
sub-tree \({\cal W}(t)
\subset {\cal W}\) defined by 
\begin{equation}\label{eqdefW(t)}
 {\cal W}(t) = \{ v \in {\cal W} : S_v < t \}.
 \end{equation}
Recalling (\ref{eqU0}) and (\ref{eqU1}),  if \(v \in \partial{\cal W}(t)\),
\[
\Upsilon_v(t-S_v)  =   \chi_0(V_v(t -S_v)) 
 + \frac{4  \tilde{m}(X_{\bar{v}})}{1-p} {\bf C}_v  \varphi(X_v,t-S_v-\tau_v){\bf 1}[\kappa_v=4,5]\cap[\tau_v \leq t-S_v].
\]
From (\ref{eqUdef2}), for \(v \in {\cal W}^\circ(t)\)
\begin{eqnarray*}
& & \Upsilon_v(t-S_v) =  \chi_0(V_v(t-S_v)) \\
& + & \frac{11m(X_{\bar{v}}) }{p} {\bf B}_v
(\Upsilon_{v \ast 0}(t -S_v- \tau_v), \Upsilon_{v \ast 1}(t -S_v - \tau_v)) {\bf 1}[ \kappa_v=1,2,3]
{\bf 1}[  \tau_v \leq t -S_v] .
\end{eqnarray*}
This is the setting of Lemma \ref{lemitcontract} with \(y_v = \chi_0(V_v(t-S_v))\), \({\bf b}_v = {\bf B}_v\), 
\(
z_v=
{\bf C}_v \varphi(X_v,t-S_v-\tau_v)\),  
\[
\eta_v = \left\{ \begin{array}{l l}
\frac{4 \tilde{m} (X_{\bar{v}})}{1-p} \quad & \mbox{if} \quad v \in \partial {\cal W}(t) \\
\\
\frac{11m(X_{\bar{v}})}{p} \quad & \mbox{if} \quad v \in {\cal W}^\circ(t) 
\end{array} \right.
\]
and 
\( \sigma_v ={\bf 1}[\kappa_v = 4,5].\)  Lemma \ref{lembivop} gives \(|{\bf B}_v(x,y)| < |x||y|\) for all \(x,y \in {\bf R}^3\).  By assumption we have \(|\chi_0(V_v(s))| \ \leq \alpha \varepsilon\) for all \(s\), and again using Lemma \ref{lemitcontract}, \(|{\bf C}_v \varphi (X_v,s )| \leq \beta \varepsilon \) for all \(s\).  The multipliers
\(\eta_v\) are bounded as follows:
\[
\frac{11m(x)}{p}= \frac{11 \int_{{\bf R}^3} |y|^{-2} h^2(x-y) dy}{8 \pi \nu p h(x)}
 \leq \frac{11 \gamma}{8 \pi \nu p} \leq \eta 
 \]
 and 
 \[
 \frac{4 \tilde{m}(x)}{1-p} = \frac{\int_{{\bf R}^3} |y|^{-1}\tilde{h}(x-y) dy}{2 \pi \nu (1-p) h(x)}
 \leq \frac{\tilde{\gamma}}{2 \pi \nu (1-p)} \leq \eta 
 \]
 for all \(x \in {\bf R}^3\).  Lemma \ref{lemitcontract} gives
  \(| \Upsilon_\phi (t)| \leq \varepsilon\) a.s. \(P_x\).  From Theorem \ref{thmbdsol} then
\[
|u(x,t)| =h(x) |E_x \Upsilon_\phi (t)| \leq h(x) E_x |\Upsilon_\phi (t)| \leq  \varepsilon h(x).
\]
\hfill \(\Box\)

\noindent{\it Proof of Theorem \ref{thmeu}:}  The existence of the weak solution \(u(x,t) = h(x) E_x \Upsilon_\phi (t)\) follows immediately from Theorem \ref{thmbdsol} and Proposition \ref{propasbd} with \(\eta =1\).  Uniqueness is derived using a martingale argument as follows.  Suppose that \(v\) is a solution to (\ref{eqexpins}) with 
\[
\sup_{x \in {\bf R}^3 , t>0} \frac{|v(x,t)|}{h(x)} \leq \varepsilon.
\]
Set \(\rho(x,t) = v(x,t)/h(x)\).  Take \(X_{\bar{\phi}}=x \in {\bf R}^3\) fixed, fix \(t>0\) and let
\begin{equation}\label{defWn(t)}
{\cal W}^{(n)}(t)=\{ v \in {\cal W}(t) : |v| \leq n\}
\end{equation}
where \({\cal W}(t)\) is as defined in (\ref{eqdefW(t)}).  For \(n \geq 0\) and \(v \in {\cal W}(t)\), 
 define the random functionals \(\Psi_v^{(n)}\) via the iterative construction
\begin{eqnarray*}
\Psi_v^{(0)}(t) & = & \chi_0(V_v(t))\\
& + & \frac{11m(X_{\bar{v}})}{p}{\bf B}_v( \rho(X_v, t-\tau_v),\rho(X_v, t-\tau_v)) 
{\bf 1}[\kappa_v =1,2,3] \cap [\tau_v \leq t]   \\
& + & \frac{4 \tilde{m}(X_{\bar{v}})}{1-p} {\bf C}_v \varphi(X_v,t-\tau_v) 
{\bf 1}[\kappa_v =4,5] \cap [\tau_v \leq t]
\end{eqnarray*}
and, for \(n \geq 1\),
\begin{eqnarray*}
\Psi_v^{(n)}(t) & = & \chi_0(V_v(t))\\
& + & \frac{11m(X_{\bar{v}})}{p}{\bf B}_v( \Psi_{v \ast 0}^{(n-1)}(t-\tau_v),\Psi_{v \ast 1}^{(n-1)}(t-\tau_v)) 
{\bf 1}[\kappa_v =1,2,3] \cap [\tau_v \leq t]   \\
& + &  \frac{4 \tilde{m}(X_{\bar{v}})}{1-p} {\bf C}_v \varphi(X_v,t-\tau_v) 
{\bf 1}[\kappa_v =4,5] \cap [\tau_v \leq t].
\end{eqnarray*}
Notice that for each \(n\), \(\Psi_\phi^{(n)}(t)\) depends only on those ensembles in \(\bf X\) indexed by \(v \in {\cal W}^{(n)}(t)\).  If \({\cal W}^{(n)}(t) = {\cal W}(t)\), then \(\Psi_\phi^{(n)}(t)=
\Upsilon_\phi(t)\).
Since \(v(x,t)\) is a solution to (\ref{eqexpins}),
\[
v(x,t)=h(x) \rho(x,t) = h(x) E_x\Psi_\phi^{(0)} (t)
\]
and for all \(v\),
\[
E_x( \Psi_v^{(0)}(t-S_v)| {\cal G}_{\bar{v}}){\bf 1}[ t-S_v \geq 0]  =
 \rho(X_{\bar{v}},t-S_v){\bf 1}[ t-S_v \geq 0] .
\] 
If for some \(n \geq 1\) and all \(v \in \cal V\), 
\(E_x(\Psi_v^{(n-1)}(t-S_v) |{\cal G}_{\bar{v}} )=
\rho(X_{\bar{v}},t-S_v)\) on the set \([ t-S_v \geq 0] \), then
\begin{eqnarray*}
& & E_x( \Psi_v ^{(n)}(t-S_v) | {\cal G}_v) {\bf 1}[ t-S_v \geq 0]  =  \chi_0(V_v (t-S_v)) \\
& + &\Big\{ (11m(X_{\bar{v}})/p){\bf B}_v( E_x(\Psi_{v \ast 0}^{(n-1)}(t-S_v - \tau_v)|{\cal G}_v ),
 E_x(\Psi_{v \ast 1}^{(n-1)}( t-S_v - \tau_v)|{\cal G}_v) ) 
{\bf 1}[\kappa_v =1,2,3]   \\
& + &  (4 \tilde{m}(X_{\bar{v}})/(1-p)) {\bf C}_v \varphi(X_v,t-S_v-\tau_v) 
{\bf 1}[\kappa_v =4,5] \Big\} {\bf 1} [\tau_v \leq t-S_v] \\
& = & \chi_0(V_v (t-S_v)) \\
& + &\Big\{ (11m(X_{\bar{v}})/p){\bf B}_v(  \rho(X_v, t-S_v - \tau_v), \rho(X_v, t-S_v - \tau_v))
{\bf 1}[\kappa_v =1,2,3]   \\
& + &  (4 \tilde{m}(X_{\bar{v}})/(1-p)) {\bf C}_v \varphi(X_v,t-S_v-\tau_v) 
{\bf 1}[\kappa_v =4,5] \Big\} {\bf 1} [\tau_v \leq t-S_v] .
\end{eqnarray*}
and, on the set \([t-S_v \geq 0]\),
\begin{eqnarray*}
& & E_x(\Psi_v^{(n)} (t-S_v)|{\cal G}_{\bar{v}})  =  E_x(E_x(\Psi_v^{(n)} (t-S_v)
|{\cal G}_v)|{\cal G}_{\bar{v}})\\
& = & \int_{{\bf R}^3} J(y,t-S_v|X_{\bar{v}}) \chi_0(y) dy \\
& + & \int_{s=0}^{t-S_v} \int_{{\bf R}^3}  \int_{{\bf R}^3} \bigg\{ m(X_{\bar{v}}) 
\bigg(f_0(s|z)f(y,z|X_{\bar{v}}) {\bf b}_1(z; \rho (X_{\bar{v}}-z,t-S_v-s), \rho(X_{\bar{v}}-z,t-S_v-s)) )\\
& + & \Big(2f_1(s|z)f(y,z|X_{\bar{v}}) -3f_1(s|y)f(y,z|X_{\bar{v}})\Big) 
{\bf b}_2(z; \rho (X_{\bar{v}}-z,t-S_v-s), \rho(X_{\bar{v}}-z,t-S_v-s) ) \bigg)\\
&  + & \tilde{m}(X_{\bar{v}}) \Big(2f_1(s|z) \tilde{f}(y,z|X_{\bar{v}}) {\bf P}_z \varphi(X_{\bar{v}}-z,t-S_v-s) \\
& & \quad \quad -  f_1(s|y) \tilde{f}(y,z|X_{\bar{v}}) (I-3 e_z e_z^t) \varphi (X_{\bar{v}}-z,t-S_v-s)\Big) \bigg\} dz dy ds \\
& = & \rho(X_{\bar{v}},t-S_v).
\end{eqnarray*}
By induction then, \(E_x(\Psi_v^{(n)}(t-S_v)|{\cal G}_{\bar{v}}){\bf 1}[t-S_v \geq 0] = \rho(X_{\bar{v}},t-S_v){\bf 1}[t-S_v \geq 0] \) for all \(n\) and \(v\).  In particular, \(E_x\Psi_\phi^{(n)}(t) = \rho(x,t)\) for all \(n\).  
 Since \({\cal W}(t)\) is finite a.s. \(P_x\), there exists with probability 1 some random \(N\) with
\( {\cal W}^{(n)}(t) = {\cal W}(t)\) and
\(
\Psi_\phi^{(n)}(t) = \Upsilon_\phi (t) \) for all \( n \geq N.
\)
The random functionals \(\Upsilon_\phi\) and \(\Psi_\phi^{(n)}\) are all bounded by \(\varepsilon\) in magnitude,
so
\[
|\rho(x,t) - \chi(x,t)|  =  |E_x\Psi_\phi^{(n)}(t) -E_x \Upsilon_\phi(t)| 
 \leq  E_x|\Psi_\phi^{(n)}(t) - \Upsilon_\phi(t)| 
 \leq  2 \varepsilon P_x({\cal W}^{(n)}(t) \neq {\cal W}(t)).
\]
This probability goes to 0 as \(n\) goes to \(\infty\). \hfill \(\Box\)

\noindent
{\bf Example 1:}   Take \(p=\alpha =1/2\) to see that if 
 \[
\sup_{x \in {\bf R}^3} |x||u_0(x)| < \nu/11 \pi^2 \quad \hbox{and} \quad 
\sup_{t>0} |g(x,t)| < \tilde{h}_0(|x|)  \; \; \hbox{for all} \; \; x \in {\bf R}^3
\]
for some \(\tilde{h}_0:[0,\infty) \to [0,\infty]\) with \(\int_{{\bf R}^3} \tilde{h}_0(|x|) dx \leq \nu^2/11 \pi \), then there exists a unique weak solution \(u(x,t)\) to (\ref{ns}) with \(\sup_{x,t} |x||u(x,t)| < 2 \nu/11 \pi^2\).

\noindent
{\bf Example 2:}  Again take \(p=\alpha =1/2\).  If there exists \(\{\mu_j:j\geq 1\}\) in \({\bf R}^3\) and \(\{p_j: j \geq 1\}\), \(p_j >0\), \(\sum p_j =1\) with \[
|u_0(x)| < ( \nu /11 \pi^2) \sum p_j|x-\mu_j|^{-1} \quad \hbox{and} \quad \sup_{t} |g(x,t)| <  \sum p_j \tilde{h}_0(|x-\mu_j|) \quad\hbox{for all} \quad x \in {\bf R}^3
\]
for some \(\tilde{h}_0\) with \(\int_{{\bf R}^3} \tilde{h}_0(|x|) dx \leq \nu^2/11 \pi  \),
then there exists a unique weak solution \(u(x,t)\) to (\ref{ns}) with
 \(\sup_{t} |u(x,t)| < (2 \nu/11 \pi^2) \sum p_j |x-\mu_j|^{-1} \) for all \(x \in{\bf R}^3\).

\noindent
{\bf Example 3:}
If \[
\sup_{x \in {\bf R}^3} (1+ |x|) |u_0(x)| < \nu/11 \pi^2 \quad \hbox{and} \quad 
\sup_{t>0} |g(x,t)| < \tilde{h}(x)  \; \; \hbox{for all} \; \; x \in {\bf R}^3
\]
for some \(\tilde{h}\) with \(0 \leq \tilde{h}(x) \leq \int_{{\bf R}^3} \tilde{h}_0(|x-y|)(2 \pi |y|)^{-1}(1+|y|)^{-3} dy\) where \(\int_{{\bf R}^3} \tilde{h}_0(|x|) dx \leq \nu^2/11 \pi \), then there exists a unique weak solution \(u(x,t)\) to (\ref{ns}) with
 \(\sup_{x,t}(1+ |x|) |u(x,t)| < 2 \nu/11 \pi^2\).

Notice that Examples 1,2 and 3 all give existence and uniqueness in a sub-space of the Marcinkiewicz space \(L^{3,\infty}\); c.f. Cannone and Karch \cite{CK03}.

\section{Non-excessive kernels and an alternate probabilistic representation. }

The following representation can be used with both non-excessive and excessive majorizing kernels. 
Let \((h,\tilde{h})\) be a majorizing kernel pair and let \({\bf X}^{(1)} =\{(X_v,Y_v,Z_v,\tau_v,\kappa_v): v \in {\cal V}\}\) be the \(\cal V\)-indexed Markov process with the transition density of \((Y_v,Z_v,\tau_v,\kappa_v)\), given that \(X_{\bar{v}} =x\), as defined in (\ref{eqjtdens})
and again let \(X_v=X_{\bar{v}}-Z_v\).  To define the stochastic recursion, let
\[
m_0(x,t) = \frac{\int_{{\bf R}^3} u_0(x-y)K(y,2 \nu t)dy}{h(x)}
\]
and let \(\chi_0\), \(\varphi\), \(m\), and \(\tilde{m}\) be as defined in equations (\ref{eqchi0def}) through (\ref{eqmtdef}).  We will also use the random operators \({\bf B}_v\) and \({\bf C}_v\) as defined in (\ref{eqBdef}) and (\ref{eqCdef}).  Now define the random recursive functional \(\Xi\) on \(\cal V\) via
\begin{eqnarray*}
\Xi_v(t)  =  m_0(X_{\bar{v}},t) & + & \frac{11m(X_{\bar{v}})}{p} {\bf B}_v(\Xi_{v \ast 0}(t-\tau_v), \Xi_{v \ast 1}(t-\tau_v)) {\bf 1}[ \kappa_v = 1,2,3]\cap[\tau_v \leq t]\\
&  + & \frac{4 \tilde{m}(X_{\bar{v}})}{1-p} {\bf C}_v \varphi(X_v, t- \tau_v){\bf 1}[ \kappa_v = 4,5]\cap[\tau_v \leq t].
\end{eqnarray*}  
If \((h,\tilde{h})\) is excessive, this random recursion essentially replaces the first term of \(\Upsilon_v\), \(\chi_0(V_v(t))\), with \(m_0(X_{\bar{v}},t)\), the conditional expectation of \(\chi_0(V_v(t))\) given \({\cal G}_{\bar{v}}\).
\begin{proposition}   Let \((h,\tilde{h})\) be a majorizing kernel pair with constant pair
\((\gamma, \tilde{\gamma}) \in (0,8 \pi \nu \eta p/11] \times [0, 2 \pi \nu \eta (1-p)]\) for some \( p \in (0,1/2]\) 
and \(\eta > 0\).
Suppose that for some \(\alpha \in [0,1)\) and \(\varepsilon, \beta \in (0, (1- \alpha)/\eta)\), 
\((u_0,g)\) is \((h,\tilde{h})\)-admissible with
\[
\sup_{x \in {\bf R}^3, t \geq 0} \frac{|\int_{{\bf R}^3} u_0(x-y)K(y,2 \nu t)dy|}{h(x)}  \leq \alpha \varepsilon
\quad
\mbox{and}
\quad
\sup_{x \in {\bf R}^3, t \geq 0} \frac{|g(x,t)|}{\tilde{h}(x)} \leq \beta \varepsilon .
\]
Then for all \(x \in {\bf R}^3\) and \(t>0\)
\[
|\Xi_\phi(t)| \leq \varepsilon \quad a.s. \, \, P_x.
\]
\end{proposition}
\noindent{\it Proof: } The proof of this proposition exactly parallels that of Proposition \ref{propasbd}.  Here the \(y_v\) of Lemma \ref{lemitcontract} is equal to \(m_0(X_{\bar{v}},t-S_v)\), so by assumption, \(|y_v |\leq \alpha \varepsilon\).  The proof is then identical.
\hfill \(\Box\)

The conditions appearing in the following corollary are somewhat simpler.
\begin{corollary}\label{cornonex}  Let \((h,\tilde{h})\) be a majorizing kernel pair with constant pair
\((\gamma, \tilde{\gamma}) \in (0,8 \pi \nu \eta p/11] \times [0, 2 \pi \nu \eta (1-p)]\) for some \( p \in (0,1/2]\) 
and \(\eta > 0\) and \[
\sup_{x \in {\bf R}^3, t>0} \frac{ \int_{{\bf R}^3} h(x-y) K(y,2 \nu t)dy}{h(x)} \leq M < \infty.
\]
Suppose that for some \(\alpha \in [0,1)\) and \(\varepsilon, \beta \in (0, (1- \alpha)/\eta)\),  \((u_0,g)\) is \((h,\tilde{h})\)-admissible with 
\[
\sup_{x \in{\bf R}^3} \frac{|u_0(x)|}{h(x)} \leq \alpha \varepsilon /M \quad \hbox{and} \quad 
\sup_{x \in{\bf R}^3, t>0} \frac {|g(x,t)|}{\tilde{h}(x)} \leq \beta \varepsilon .
\]
Then for all \(x \in {\bf R}^3\) and \(t>0\)
\[
|\Xi_\phi(t)| \leq \varepsilon \quad a.s. \, \, P_x.
\]
\end{corollary}

The proofs of the following theorems are omitted due to their similarity to proofs appearing in the previous section.  
\begin{theorem}\label{thmbdnonex}  Let \((h, \tilde{h})\) be a majorizing kernel pair.  
If  
\(
E_x| \Xi_\phi(t)| \leq M < \infty \) for all \(x\) and \(t\), then 
\[ u(x,t) = h(x) E_x \Xi_\phi(t)\]
is a weak solution to the Navier-Stokes equations.
\end{theorem}

\begin{theorem}\label{thmgeneu}  (Existence and Uniqueness)  Let \((h,\tilde{h})\) be a majorizing kernel pair with constant pair \((\gamma, \tilde{\gamma}) \in (0, 8 \pi \nu p/11] \times [0, 2 \pi \nu (1-p)]\) for some \(p \in (0,1/2]\).  If \((u_0,g)\) is \((h,\tilde{h})\)-admissible with
\begin{equation}\label{eqnecond}
\sup_{x \in {\bf R}^3, t \geq 0} \frac{|\int_{{\bf R}^3} u_0(x-y)K(y,2 \nu t)dy|}{h(x)}  \leq \alpha \varepsilon
\quad
\mbox{and}
\quad
\sup_{x \in {\bf R}^3, t \geq 0} \frac{|g(x,t)|}{\tilde{h}(x)} \leq \beta \varepsilon .
\end{equation}
for some \(\alpha \in [0,1)\) and some \(\varepsilon, \; \beta \in (0,1- \alpha)\),
then 
\(
u(x,t)= h(x) E_x \Xi_\phi (t)
\)
is a weak solution to the Navier-Stokes equations with
\[
\sup_{x \in {\bf R}^3, t \geq 0} \frac{|u(x,t)|}{h(x)} \leq \varepsilon.
\]
This solution is unique in the class 
\(
\{ v \in ({\cal S}' ({\bf R}^3 \times (0,\infty)))^3 :
 \sup_{x \in {\bf R}^3, t \geq 0} \frac{|v(x,t)|}{h(x)} \leq \varepsilon\}.
\)
\end{theorem}

Three useful corollaries follow.  
The first imposes separate bounds on the magnitude of the initial velocity relative to \(h\) and the size of \(h*K\) relative to \(h\).   In the last two corollaries a single bound is imposed on the magnitude of \(u_0*h\) relative to \(h\).  The last corollary treats the case of no forcing.  This is of particular interest in regards to regularity; see remark (4) in the following section.
\begin{corollary}\label{corneeu}  Let \((h,\tilde{h})\) be a majorizing kernel pair with 
\begin{equation}\label{eqcondM}
 \sup_{x \in {\bf R}^3, t>0} \frac{ \int_{{\bf R}^3} h(x-y) K(y,2 \nu t)dy}{h(x)} \leq M < \infty.
\end{equation}
and constant pair \((\gamma, \tilde{\gamma}) \in (0, 8 \pi \nu p/11] \times [0, 2 \pi \nu (1-p)]\) for some \(p \in (0,1/2]\).
 If \((u_0,g)\) is \((h,\tilde{h})\)-admissible with
\begin{equation}\label{eqMratio}
\sup_{x \in {\bf R}^3} \frac{|u_0(x)|}{h(x)}  \leq \alpha \varepsilon/M
\quad
\mbox{and}
\quad
\sup_{x \in {\bf R}^3, t \geq 0} \frac{|g(x,t)|}{\tilde{h}(x)} \leq \beta \varepsilon 
\end{equation}
for some \(\alpha \in [0,1)\) and \(\varepsilon, \; \beta \in (0,1- \alpha)\).
Then 
\(
u(x,t)= h(x) E_x \Xi_\phi (t)
\)
is a weak solution to the Navier-Stokes equations with
\[
\sup_{x \in {\bf R}^3, t \geq 0} \frac{|u(x,t)|}{h(x)} \leq \varepsilon.
\]
This solution is unique in the class 
\(
\{ v \in ({\cal S}' ({\bf R}^3 \times (0,\infty)))^3 :
 \sup_{x \in {\bf R}^3, t \geq 0} \frac{|v(x,t)|}{h(x)} \leq \varepsilon\}.
\)
\end{corollary}
\noindent{\it Proof of Corollary \ref{corneeu}:}  Combining (\ref{eqcondM}) and (\ref{eqMratio}) gives (\ref{eqnecond}). \hfill \(\Box\)

The next corollary relies on the rescaling seen in the derivation of Corollary \ref{coreu}.
\begin{corollary}\label{corprethm1}  Let \((h,\tilde{h})\) be a majorizing kernel pair with constant pair \((\gamma, \tilde{\gamma})\), \(\tilde{\gamma} >0\), and suppose that \((u_0,g)\) is \((h,\tilde{h})\)-admissible with
\begin{equation}\label{eqcorprethm1}
\sup_{x \in {\bf R}^3, t \geq 0} \frac{|\int_{{\bf R}^3} u_0(x-y)K(y,2 \nu t)dy|}{h(x)}  \leq 
\frac{8 \pi \nu p \alpha \varepsilon}{11 \gamma}
\quad
\mbox{and}
\quad
\sup_{x \in {\bf R}^3, t \geq 0} \frac{|g(x,t)|}{\tilde{h}(x)} \leq 
\frac{(4 \pi \nu)^2 p (1-p) \beta \varepsilon}{11 \gamma \tilde{\gamma}} 
\end{equation}
for some \(p \in [0,1/2]\), \(\alpha \in [0,1)\) and \(\varepsilon, \; \beta \in (0,1- \alpha)\),
then 
\(
u(x,t)= h(x) E_x \Xi_\phi (t)
\)
is a weak solution to the Navier-Stokes equations with
\begin{equation}\label{eqresult}
\sup_{x \in {\bf R}^3, t \geq 0} \frac{|u(x,t)|}{h(x)} \leq \frac{8 \pi \nu p \varepsilon}{11 \gamma}.
\end{equation}
This solution is unique in the class 
\(
\{ v \in ({\cal S}' ({\bf R}^3 \times (0,\infty)))^3 :
 \sup_{x \in {\bf R}^3, t \geq 0} \frac{|v(x,t)|}{h(x)} \leq \frac{8 \pi \nu p \varepsilon}{11 \gamma}\}.
\)
\end{corollary}
\noindent{\it Proof of Theorem \ref{thm1}:}  Equation (\ref{eqconth1}) guarantees that (\ref{eqcorprethm1}) is satisfied with \(\gamma=\tilde{\gamma}=1\) and \(p=\alpha =1/2\).
Then (\ref{eqresult}) holds for any \(\varepsilon < 1/2\).  \hfill \(\Box\)

\begin{corollary}  Let \((h,0)\) be a majorizing kernel pair with constant pair \((\gamma, 0)\) and suppose that 
\[
\sup_{x \in {\bf R}^3, t \geq 0} \frac{|\int_{{\bf R}^3} u_0(x-y)K(y,2 \nu t)dy|}{h(x)}  \leq 
\frac{4 \pi \nu \alpha \varepsilon}{11 \gamma}
\]
for some \(\alpha \in [0,1)\) and \(\varepsilon  \in (0,1- \alpha)\).
Then 
\(
u(x,t)= h(x) E_x \Xi_\phi (t)
\)
is a weak solution to the Navier-Stokes equations with
\begin{equation}
\sup_{x \in {\bf R}^3, t \geq 0} \frac{|u(x,t)|}{h(x)} \leq \frac{4 \pi \nu \varepsilon}{11 \gamma}.
\end{equation}
This solution is unique in the class 
\(
\{ v \in ({\cal S}' ({\bf R}^3 \times (0,\infty)))^3 :
 \sup_{x \in {\bf R}^3, t \geq 0} \frac{|v(x,t)|}{h(x)} \leq \frac{4 \pi \nu \varepsilon}{11 \gamma}\}.
\)
\end{corollary}

\noindent
{\bf Example 4:}   If \(\sup_{x} (1+|x|^2)|u_0(x)| \leq 11 \pi/16 \nu (1+3e^{-2/3})\) and the forcing \(g\) is identically 0, then there exists a unique weak solution to (\ref{ns}) with 
\(\sup_{x,t}(1+|x|^2)|u(x,t)| < 11 \pi / 8 \nu\).

\noindent
{\bf Example 5:}  Taking \(H_p(x)=(1+|x|)^{-p}\) for \(p \in (1,2]\) fixed, \[
\int_{{\bf R}^3} H_p(x-y)K(y,2 \nu t)dy \leq (1+3e^{-2/3})^{p/2}(1+1/\sqrt{2})^{p/2} H_p(x) 
\]
for all \(x \in {\bf R}^3\) and \(t >0\).  Then, if 
\(\sup_{x} (1+|x|)^p |u_0(x)| \leq  \pi^{p-3} \nu (1+3e^{-2/3})^{-p/2} (1+1/\sqrt{2})^{1-3p/2}/11\)
and the forcing \(g\) is identically 0, there exists a unique weak solution to (\ref{ns}) with 
\(\sup_{x,t}(1+|x|)^p|u(x,t)| < 2 \pi^{p-3} \nu (1+1/\sqrt{2})^{1-p}/11\).

\section{Remarks.} 

(1)  It appears that the integral formulation on which
 the probabilistic representation is based can also be used as the basis for 
Picard iteration.  The advantage of the probabilistic approach as it appears here is that it allows, with probability 1, the assumption that there are a finite number of iterates.  This permits easy derivation of bounds via inductive arguments.
 
 (2)  The probabilistic representations of the Fourier transformed Navier-Stokes equations given by
Bhattacharya et al \cite{BCDGOOTW03} give both short time existence and uniqueness of solutions for large initial data and decay over time of the Fourier transformed solution.  It is possible that a physical space representation related to the ones given in this paper will also indicate short time existence and uniqueness as well as the time decay of solutions. 
 
 (3)  A semi-group approach underlies this probabilistic representation.  Montgomery-Smith \cite{M-S01}
 uses such an approach to demonstrate blow-up in finite time of an equation related to the Navier-Stokes equations.  A more complete understanding of the integration of the bilinear forms \({\bf b}_1\) and \({\bf b}_2\) may aid in understanding possible blow-up of the Navier-Stokes equations in \({\bf R}^3\). 

(4)  In some cases the results given in this paper can be used to demonstrate regularity and consequently uniqueness in a larger class of functions.  The argument is as follows: assume that \(u_0 \in L_2^3({\bf R}^3)\) is regular and take the forcing to be identically 0.  Take \((h,0)\) to be a majorizing kernel pair with 
\(h \in L_2({\bf R}^3) \cap L_q({\bf R}^3)\) for some \(q >3\).  If 
\(\sup_{x,t} | \int u_0(x-y)K(y,2 \nu t) dy |/h(x)\) is small enough, then \( \sup_{x,t}|u(x,t)|/h(x) < \varepsilon\) for suitable \(\varepsilon >0\) and for any \(r>0\) and \(T< \infty\)
\[
\int_0^T (\int|u(x,t)|^qdx)^r dt < \varepsilon^{qr} \int_0^T(\int h^q(x)dx)^rdt < \infty .
\]
That is, the Ladyzhenskaya-Prodi-Serrin condition is satisfied, and the solution \(u\) is regular and unique in a larger function space; c.f. Temam \cite{T79}, Chapter III.  A suitable kernel here is \(h=H_p\) of Example 5 with \(p \in (3/2,2]\).

 \end{document}